\documentclass[sort&compress,3p,fleqn]{elsarticle}
\usepackage{pdfpages,color,tikz,pdfsync,enumitem}
\usepackage{graphicx}
\usepackage{amsmath,amssymb,amsthm}
\usepackage[CJKbookmarks=True]{hyperref}

\allowdisplaybreaks[4]

\numberwithin{equation}{section}

\newtheorem{Def}{Definition}[section]
\newtheorem{Thm}[Def]{Theorem}
\newtheorem{Lem}[Def]{Lemma}
\newtheorem{Rmk}[Def]{Remark}

\newcommand{\R}{\mathbb{R}}
\newcommand{\J}{\mathbb{J}}
\newcommand{\mC}{\mathbb{C}}
\newcommand{\Z}{\mathbb{Z}}

\newcommand{\I}{i}

\newcommand{\ud}{d}

\newcommand{\re}{\mathrm{Re}\,}
\newcommand{\im}{\mathrm{Im}\,}

\newcommand{\Hess}{\mathrm{Hess}}

\newcommand{\bvec}[1]{\boldsymbol{#1}}

\newcommand{\vj}{\bvec{j}}

\newcommand{\vn}{\bvec{n}}
\newcommand{\vx}{\bvec{x}}
\newcommand{\va}{\bvec{a}}
\newcommand{\ta}{\widetilde{\bvec{a}}}

\newcommand{\vb}{\bvec{b}}
\newcommand{\tb}{\widetilde{\bvec{b}}}

\newcommand{\vt}{\bvec{\tau}}
\newcommand{\p}{\partial}

\newcommand{\la}{\left\langle}
\newcommand{\ra}{\right\rangle}
\newcommand{\nn}{\nonumber\\}

\newcommand{\SE}{Schr\"{o}dinger equation}

\newcommand{\wto}{\rightharpoonup}

\newcommand{\be}{\begin{equation}}
\newcommand{\ee}{\end{equation}}
\begin{document}

\begin{frontmatter}
\title{Quantized Vortex Dynamics of  the Coupled  Nonlinear Schr\"{o}dinger Equation}

\author[Zhu]{Yongxing Zhu\texorpdfstring{\corref{cor}}{}}
\ead{zhuyongxing@hnu.edu.cn}

\address[Zhu]{School of Mathematics, Hunan University, Changsha City, Hunan Province, 410082, China}
\cortext[cor]{Corresponding author.}

\begin{abstract}
% This is the abstract.
We derive rigorously the reduced dynamical law for quantized vortex dynamics of the coupled nonlinear \SE{} without Josephson junction (CNLS) when the core size of vortex $\varepsilon\to 0$.  It is proved that when $\varepsilon\to 0$, the vortex motion of one component won't affect the vortex motion on the other component. Moreover, the motion of vortices of each component follows the vortex motion law for the nonlinear \SE. 
\end{abstract}

\begin{keyword}
coupled nonlinear \SE{} \sep quantized vortex \sep canonical harmonic map \sep reduced dynamical law \sep renormalized energy \sep vortex path
% \MSC[2020]{35Q40, 35Q55}
\end{keyword}

% \date{}
\end{frontmatter}
\section{Introduction}\label{sec:Introduction}
% This is a citation\cite{AguarelesChapmanWitelski2008INteractionSprialCGL,AguarelesChapmanWitelski2010MotionSpiralCGL}.
In this article, we consider the reduced dynamical law for the vortices of the coupled \SE{} without Josephson junction:
\begin{equation}\label{eq:main}
\begin{cases}
\I \p_tu^\varepsilon-\Delta u^\varepsilon+\frac{1}{\varepsilon^2}(|u^\varepsilon|^2+g|v^\varepsilon|^2-1)u^\varepsilon=0,&\vx\in\Omega,t>0,\\
\I \p_tv^\varepsilon-\Delta v^\varepsilon+\frac{1}{\varepsilon^2}(g|u^\varepsilon|^2+|v^\varepsilon|^2-1)v^\varepsilon=0,&\vx\in\Omega,t>0,\\
u^\varepsilon(\vx,0)=u^\varepsilon_0(\vx),\quad v^\varepsilon(\vx,0)=v^\varepsilon(\vx),&\vx\in\Omega,\\
\frac{\p u}{\p \bvec{n}}(\vx)=\frac{\p v}{\p \bvec{n}}(\vx)=0,&\vx\in\p\Omega.
\end{cases}\tag{CNLS}
\end{equation}
\addtocounter{equation}{1}
Here, $u^\varepsilon,v^\varepsilon,u^\varepsilon_0,v^\varepsilon_0$ are complex-valued functions, $\Omega\subset\R^2$ is a bounded domain, $\vx=(x,y)^T$ is the spatial coordinate, $t$ is the time variable, $0<\varepsilon<<1$ is a constant which characterizes the vortex core size, and $0<g<1$ is a real parameter.
 % For simplicity, we assume $g_{11}=g_{22}=1,g_{12}=g,\omega=0$.

We define the energy for \eqref{eq:main} by
\begin{equation}\label{eq:def of E}
	E^\varepsilon_0(u^\varepsilon,v^\varepsilon):=\int_\Omega e^\varepsilon_0(u^\varepsilon,v^\varepsilon)\ud\vx,
\end{equation}
and the momentum by
\begin{equation}\label{eq:def of Q}
	\bvec{Q}(u^\varepsilon,v^\varepsilon):=\int_\Omega (\vj(u^\varepsilon)+\vj(v^\varepsilon))\ud\vx.
\end{equation}
Here, the energy density $e^\varepsilon_0$ and current $\vj$ is defined as follows: for any complex-valued functions $u,v:\Omega\to\mC$,
\begin{equation}\label{eq:def of e0e}
	e^\varepsilon_0(u,v):=\frac{1}{2}|\nabla u|^2+\frac{1}{2}|\nabla v|^2+\frac{1}{4\varepsilon^2}\left(|u|^2-\frac{1}{1+g}\right)^2+\frac{1}{4\varepsilon^2}\left(|v|^2-\frac{1}{1+g}\right)^2+\frac{g}{2\varepsilon^2}\left(|u|^2-\frac{1}{1+g}\right)\left(|v|^2-\frac{1}{1+g}\right)
\end{equation}
and 
\begin{equation}
	\vj(u)=(j_1(u),j_2(u))^T:=\im(\overline{u}\nabla u).
\end{equation}
It's well-known that \eqref{eq:main} conserves both the energy and the momentum.

{%\color{red}
The coupled \SE{} is widely used to describe the multi-component Bose-Einstein condensation \cite{BaoCai2018SpinorBEC}. Bose-Einstein condensate (BEC) is predicted by Einstein in 1924 and realized via experiment in 1995 \cite{Fadel2021MultiComponentBEC}. The early experiments focused on the case where the condensed bosons don't have any internal degree of freedom and these phenomena can be described by the nonlinear \SE:
\begin{equation}\label{eq:GPE}
	i\p_t \varphi^\varepsilon-\Delta \varphi^\varepsilon+\frac{1}{\varepsilon^2}(|\varphi^\varepsilon|^2-1)\varphi^\varepsilon=0.
\end{equation}
In recent experiments, BEC where the particles in different internal states can be distinguished, i.e. the multi-component BEC, is realized \cite{Fadel2021MultiComponentBEC}. To describe the multi-component BEC, the nonlinear \SE{} is generalized to various multi-component nonlinear \SE{} (e.g. \eqref{eq:main}, spin-F BEC).

In the study of BEC and other phenomena related to BEC (e.g. superfluidity and superconductivity), quantized vortex is a key signature.  In two dimensions, quantized vortices are particle-like defects whose centers are zeros of the order parameter. It was proved that the vortex is dynamically stable only if the winding number of the vortex is $1$ of $-1$. When considering \eqref{eq:main}, there are some differences on the definition of vortex and the winding number of the vortex. A vortex of \eqref{eq:main} should be a zero point of $u^\varepsilon$ or $v^\varepsilon$. Formally, when $\varepsilon<<1$ and $\vx$ is away from the vortices, we have 
\[
	|u^\varepsilon(\vx)|\approx\frac{1}{\sqrt{1+g}}\approx|v^\varepsilon(\vx)|.
\]
Hence, we can formally obtain 
\begin{equation}
	(u^\varepsilon(\vx),v^\varepsilon(\vx))\approx\frac{1}{\sqrt{1+g}}e^{i\theta^\varepsilon_1(\vx)}(e^{i\theta^\varepsilon_2(\vx)},e^{-i\theta^\varepsilon_2(\vx)}).
\end{equation}
If we focus on the different part of the phases of $u^\varepsilon,v^\varepsilon$, i.e. $\theta_2^\varepsilon$, we will find its winding number around a vortex would be a half integer. Hence the vortex appearing in \eqref{eq:main} is also called the fractional vortex.

The dynamics of vortices of \eqref{eq:GPE} has been widely studied. We refer to \cite{CollianderJerrard1999GLvorticesSE,Lin1998GinzburgLandauDyanmicsVorticesFilamentsSubmanifolds,ZhuBaoJian2023,BaoTang2014NumericalStudyQuantizedVorticesSchrodinger} and the references therein for the results on bounded domains, the whole plane $\R^2$ and the unit torus.  An interesting fact is that we can describe the motion of vortices via an ODE. More precisely speaking, if the initial data of \eqref{eq:GPE} possesses vortices $\va^{\varepsilon,0}_1,\cdots,\va^{\varepsilon,0}_M$  with winding number $d_1,\cdots,d_M$ and  $\va_j^0:=\lim_{\varepsilon\to 0}\va_j^{\varepsilon,0}$ exist. Then the solution $\varphi^\varepsilon(\vx,t)$ also possess vortices $\va^\varepsilon_j(t),j=1,\cdots,M$ with winding numbers $d_j,j=1,\cdots,d_M$, and the limits $\va_j(t):=\lim_{\varepsilon\to 0}\va^\varepsilon_j(t)$ exist. Meanwhile, $\va=\va(t)=(\va_1(t),\cdots,\va_M(t))$ satisfies the following ODE:
\begin{equation}
	\dot{\va}_j=-\frac{d_j}{\pi}\J\nabla _{\va_j}W_{d}(\va),
\end{equation}
where $W_{d}$ is the renormalized energy and 
\begin{equation}
	\J:=\left(\begin{array}{cc}0&1\\-1&0 \end{array}\right).
\end{equation}
The specific definition of the renormalized energy is related to the boundary condition and the winding number $d=(d_1,\cdots,d_M)$, and it usually take the following form:
\begin{equation}
	W_d(\va)=-\pi\sum_{1\le j\ne k\le M}d_jd_k\log|\va_j-\va_k|+\text{term related to the boundary condition}
\end{equation}
 As for the vortices for \eqref{eq:main}, there are much less results. Lin-Lin \cite{LinLin2003CoupledGinzburgLandau}, Lin-Lin-Wei \cite{LinLinWei2009SkyrmionsGrossPitaevskii} studied special cases where there are only two vortices. Eto-Ikeno-Nitta\cite{EtoIkenoNitta2020Fractionalvortices}, Bao-Cai\cite{BaoCai2018SpinorBEC} gave some numerical results on the dynamical laws of vortex. But to the knowledge of the author, there are still no article which give a specific ODE describing the motion of fractional vortices. 
In this article, we will give a theoretical result on the dynamical law of the vortices for the solutions of \eqref{eq:main}, i.e. we will give an ODE which describe the limit behavior of vortices of \eqref{eq:main}, just like the results for \eqref{eq:GPE} \cite{CollianderJerrard1999GLvorticesSE,LinXin1999IncompressibleLimitVortexSchrodinger,BaoTang2014NumericalStudyQuantizedVorticesSchrodinger}. 
}

To state our main result, we have to introduce some notations.
We define
\begin{equation}\label{eq:def of gammag}
	\gamma_g:=\lim_{\varepsilon\to 0^+}\left(\inf_{(u,v)\in A_g(B_1(\bvec{0}))}\int_{B_1(\bvec{0})}e^\varepsilon_0(u,v)\ud\vx-\frac{\pi}{1+g}\log\frac{1}{\varepsilon\sqrt{1+g}}\right),
\end{equation}
with
\[
	A_g(B_1(\bvec{0}))=\left\{(u,v)\in H^1(B_1(\bvec{0}))\times H^1(B_1(\bvec{0}))\left|(u(\vx),v(\vx))=\frac{1}{\sqrt{1+g}}\left(\frac{z}{|z|},1\right)\  \text{for}\ \vx\in \p B_1(\bvec{0})\right.\right\}.
\]
The proof of the existence of $\gamma_g$ is similar to the existence of $\gamma$ defined in \cite{JerrardSpirn2008RefinedEstimateGrossPitaevskiiVortices}. To make this article complete, we will give the proof in Section \ref{sec:ODE}.

We denote
\[
	\Omega^M_*:=\{(\vx_1,\cdots,\vx_M)\in \Omega^M|\vx_j\ne\vx_k\ \text{for any}\ j\ne k\}.
\]
Then for $\va=(\va_1,\cdots,\va_M)\in \Omega^M_*$ and $d=(d_1,\cdots,d_M)\in\{\pm 1\}^M$, we define the renormalized energy $W_d(\va)$ by
\begin{equation}
	W_d(\va)=-\pi\left(\sum_{1\le j\ne k\le M}d_jd_k\log|\va_j-\va_k|+\sum_{j,k=1}^M d_jd_kF(\va_j,\va_k)\right),
\end{equation}
where $F(\vx,\bvec{y})$ solves
\begin{equation}
	\Delta_{\vx}F(\vx,\bvec{y})=0\ \text{in}\ \Omega,\quad F(\vx,\bvec{y})=-\log|\vx-\bvec{y}|\ \text{for any}\  \vx\in\p\Omega,\bvec{y}\in\Omega.
\end{equation}
And we define 
\begin{equation}
	J(u):=\im(\p_xu\p_y\overline{u})=\frac{1}{2}\nabla\cdot(\J\vj(u)).
\end{equation}
Then, we state our main result as follows:
\begin{Thm}\label{thm:main}
Assume that there exist $(\va^0_1,\cdots,\va^0_M,\vb^0_1,\cdots,\vb_N^0)\in \Omega^{M+N}_*,$ $d^1=(d_1^1,\cdots,d_M^1)\in\{1,-1\}^M ,d^2=(d_1^2,\cdots,d_N^2)\in\{1,-1\}^N$, such that 
\begin{equation}\label{eq:convergence of Ju0}
	J(u_0^\varepsilon)\to \frac{\pi}{1+g}\sum_{j=1}^Md_j^1\delta_{\va_j^0},\quad J(v_0^\varepsilon)\to \frac{\pi}{1+g}\sum_{j=1}^Nd_j^2\delta_{\vb_j^0},\quad \text{in}\  W^{-1,1}(\Omega):=(C^1(\Omega))'
\end{equation}
\begin{equation}\label{eq:energy bound of initial data}
E^\varepsilon_0(u^\varepsilon_0,v^\varepsilon_0)=(M+N)\left(\frac{\pi}{{1+g}}\log\frac{1}{\sqrt{1+g}\varepsilon}+\gamma_g\right)+\frac{1}{1+g}W_{d^1}(\va^0)+\frac{1}{1+g}W_{d^2}(\vb^0)+o(1),
\end{equation}
then there exist $T>0$ and Lipschitz paths $\va_j,\vb_k:[0,T)\to \Omega,1\le j\le M,1\le k\le N$ such that 
\begin{equation}\label{eq:convergence of Ju(t)}
	J(u^\varepsilon(t))\to \frac{\pi}{1+g}\sum_{j=1}^Md_j^1\delta_{\va_j(t)},\quad J(v^\varepsilon)\to \frac{\pi}{1+g}\sum_{j=1}^Nd_j^2\delta_{\vb_j(t)},\quad \text{in}\  W^{-1,1}(\Omega)
\end{equation}
and $\va=\va(t)=(\va_1(t),\cdots,\va_M(t)),\vb=\vb(t)=(\vb_1(t),\cdots,\vb_N(t))$ satisfy
\begin{equation}\label{eq:MainODE}
\begin{cases}
\dot{\va}_j=-\frac{d_j^1}{\pi}\J\nabla_{\va_j}W_{d^1}(\va),\\
\dot{\vb}_j=-\frac{d_j^2}{\pi}\J\nabla_{\vb_j}W_{d^2}(\vb),\\
\va(0)=\va^0=(\va^0_1,\cdots,\va_M^0),\quad\vb(0)=\vb^0=(\vb_1^0,\cdots,\vb_N^0).
\end{cases}
\end{equation}
\end{Thm} 

The rest of this paper is organized as follows. In Section \ref{sec:preliminaries}, we will notations and estimates as a preparation for the proof of the main theorem. In Section \ref{sec:proof of main theorem}, we will prove the main theorem. We will first prove the existence of vortices, then prove the dynamical law. In Section \ref{sec:lower bound near vortex}, we will study the energy near a vortex. In Section \ref{sec:ODE}, we will prove the existence of $\gamma_g$.
\section{Notations and estimates}\label{sec:preliminaries}

For two complex numbers $z,w$, we define 
\[
	\la z,w\ra:=\frac{\overline{z}w+\overline{w}z}{2}.
\]
More generally, we define 
\begin{equation}
	\la\bvec{z},\bvec{w}\ra:=\sum_{j=1}^M\la z_j,w_j\ra
\end{equation}
for vectors $\bvec{z}=(z_1,\cdots,z_M)^T,\bvec{w}=(w_1,\cdots,w_M)^T\in \mC^M$.

For any $\va=(\va_1,\cdots,\va_M)\in \Omega^M,\vb=(\vb_1,\cdots,\vb_N)\in \Omega^N$, we define
\begin{equation}
	r(\va,\vb)=\frac{1}{4}\min\{r_{\va},r_{\vb},r_{\va\vb}\}
\end{equation}
with 
\[
	r_{\va}=\min\left\{\min_{1\le j<k\le M}|\va_j-\va_k|,\min_{1\le j\le M}dist(\va_j,\p\Omega)\right\},
\]
\[
	r_{\vb}=\min\left\{\min_{1\le j<k\le N}|\vb_j-\vb_k|,\min_{1\le j\le N}dist(\vb_j,\p\Omega)\right\},
\]
and
\[
	r_{\va\vb}=\min_{1\le j\le M,1\le k\le N}|\va_j-\vb_k|.
\]
If $r(\va,\vb)>0$, we define
\begin{equation}
	\Omega_\rho(\va,\vb):=\Omega\setminus\left(\left(\cup_{j=1}^MB_\rho(\va_j)\right)\cup\left(\cup_{j=1}^NB_\rho(\vb_j)\right) \right)
\end{equation}
for any $\rho<r(\va,\vb)$. Then we define
\begin{equation}
	\Omega_*(\va,\vb):=\cup_{\rho>0}\Omega_\rho(\va,\vb)=\Omega\setminus\{\va_1,\cdots,\va_M,\vb_1\cdots,\vb_N\}.
\end{equation}

Following \cite{CollianderJerrard1999GLvorticesSE}, we define 
\begin{equation}
	E^\varepsilon(u):=\int_\Omega e^\varepsilon(u)\ud\vx,\quad e^\varepsilon(u)=\frac{1}{2}|\nabla u|^2+\frac{1}{4\varepsilon^2}(|u|^2-1)^2.
\end{equation}
Via the definition \eqref{eq:def of E}, there is a simple estimate:
\begin{align}\label{eq:Egu+Egv<Euv}
&E^{g_-\varepsilon}(\sqrt{1+g}u)+E^{g_-\varepsilon}(\sqrt{1+g}v)\nn \le &E^{\varepsilon\sqrt{1+g}}(\sqrt{1+g} u)+E^{\varepsilon\sqrt{1+g}}(\sqrt{1+g} v)+\frac{g}{2\varepsilon^2}\int_\Omega\left(|u|^2-\frac{1}{1+g}\right)\left(|v|^2-\frac{1}{1+g}\right)\ud\vx\nn 
=&(1+g)E^\varepsilon_0(u,v)\le E^{\varepsilon}(\sqrt{1+g}u)+E^{\varepsilon}(\sqrt{1+g}v),
\end{align}
with $g_-=\sqrt{(1-g)/(1+g)}$.
% , and 
% \begin{equation}
% 	E^\varepsilon(u)=\int_\Omega e^\varepsilon(u)\ud\vx,\quad e^\varepsilon(u)=\frac{1}{2}|\nabla u|^2+\frac{1}{4\varepsilon^2}(|u|^2-1)^2.
% \end{equation}

We define the canonical harmonic map $H$. Following \cite{JerrardSpirn2008RefinedEstimateGrossPitaevskiiVortices}, for $\va=(\va_1,\cdots,\va_M)\in\Omega^M_*$ and $d=(d_1,\cdots,d_M)\in\{\pm1\}^M$, the canonical map $H=H(\vx)=H_d(\vx;\va)$ is the solution of 
\begin{equation}\label{eq:def of H initial}
	\nabla \cdot \vj(H)=0,\nabla\cdot(\J\vj(H))=2\pi\sum_{j=1}^Md_j\delta_{\va_j},
\end{equation}
with boundary condition $\vn\cdot\vj(H)=\bvec{0}$. Then \cite{JerrardSpirn2008RefinedEstimateGrossPitaevskiiVortices} gives
\begin{equation}\label{eq:energy of eH}
	\int_{\Omega_\rho(\va)}|\nabla H|^2\ud\vx-2M\pi\log\frac{1}{\rho}-W_d(\va)=O(\rho^2).
\end{equation}

For $\va=(\va_1,\cdots,\va_M)\in \Omega^M,\vb=(\vb_1,\cdots,\vb_N)\in\Omega^N$ such that $r(\va,\vb)>0$ and $d^1=(d^1_1,\cdots,d^1_M),d^2=(d_1^2,\cdots,d_N^2)$,we can define $u_*=\frac{1}{\sqrt{1+g}}H_{d^1}(\vx;\va),v_*=\frac{1}{\sqrt{1+g}}H_{d^2}(\vx;\vb)$. Then \eqref{eq:def of e0e} and \eqref{eq:energy of eH} imply
\begin{equation}\label{eq:estimate of euv*}
	\int_{\Omega_\rho(\va,\vb)}e^\varepsilon_0(u_*,v_*)=\frac{1}{1+g}\left((M+N)\pi\log\frac{1}{\rho}+W_{d^1}(\va)+W_{d^2}(\vb) \right)+O(\rho^2).
\end{equation}

Noting that \begin{equation}
	\vj(u_*(\vx))=\frac{1}{1+g}\vj(H_{d^1}(\vx;\va)),\vj(v_*(\vx))=\frac{1}{1+g}\vj(H_{d^2}(\vx;\va)),
\end{equation}
Lemma 2.3.1 in \cite{CollianderJerrard1999GLvorticesSE} implies
\begin{equation}\label{eq:hessphijujv}
	\int_\Omega (\la\Hess (\varphi)\vj(u^*),\J\vj(u^*)\ra+\la\Hess (\varphi)\vj(v^*),\J\vj(v^*)\ra)\ud\vx=-\frac{1}{(1+g)^2}\nabla\varphi(\va_j)\cdot(\J\nabla_{\va_j}W_{d^1}(\va)),
\end{equation}
where $\varphi\in C_0^\infty(B_{r(\va,\vb)}(\va_j))$ and there is $\bvec{p}=(p_1,p_2)^2\in\R^2$ such that
\begin{equation}
	\nabla \varphi(\vx)=\bvec{p},\vx\in B_R(\va_j) \ \text{for some}\ 0<R<r(\va,\vb).
\end{equation}
\section{The proof of the main theorem}\label{sec:proof of main theorem}
\subsection{Derivative of quantities related to the solution}
Taking derivative of $e^\varepsilon_0(u^\varepsilon,v^\varepsilon),\vj(u^\varepsilon),\vj(v^\varepsilon),|u^\varepsilon|^2,|v^\varepsilon|^2$ and noting \eqref{eq:def of E}, \eqref{eq:main}, we have
\begin{equation}
	\frac{\ud}{\ud t}e^\varepsilon_0(u^\varepsilon,v^\varepsilon)=\nabla\cdot\re(\overline{u_t^\varepsilon}\nabla u^\varepsilon)+\nabla\cdot\re(\overline{v_t^\varepsilon}\nabla v^\varepsilon),
\end{equation}
which immediately implies
\begin{equation}\label{eq:conserve energy}
	E^\varepsilon_0(u^\varepsilon,v^\varepsilon)\equiv E^\varepsilon_0(u_0^\varepsilon,v_0^\varepsilon),
\end{equation}
\begin{equation}\label{eq:conserve mass}
	\frac{\ud}{\ud t}|u^\varepsilon|^2=2\nabla\cdot\vj(u^\varepsilon),\quad \frac{\ud}{\ud t}|v^\varepsilon|^2=2\nabla\cdot\vj(v^\varepsilon),
\end{equation}
\begin{equation}
	\frac{\ud}{\ud t}\vj(u^\varepsilon)=2\nabla\cdot(\nabla u^\varepsilon\otimes\nabla u^\varepsilon)+\nabla\left(\im(\overline{u^\varepsilon}u_t^\varepsilon)-|\nabla u^\varepsilon|^2-\frac{1}{2\varepsilon^2}|u^\varepsilon|^4-\frac{1}{\varepsilon^2}|u^\varepsilon|^2 \right)-\frac{g}{\varepsilon^2}|v^\varepsilon|^2\nabla|u^\varepsilon|^2,
\end{equation}
\begin{equation}
	\frac{\ud}{\ud t}\vj(v^\varepsilon)=2\nabla\cdot(\nabla v^\varepsilon\otimes\nabla v^\varepsilon)+\nabla\left(\im(\overline{v^\varepsilon}v_t^\varepsilon)-|\nabla v^\varepsilon|^2-\frac{1}{2\varepsilon^2}|v^\varepsilon|^4-\frac{1}{\varepsilon^2}|v^\varepsilon|^2 \right)-\frac{g}{\varepsilon^2}|u^\varepsilon|^2\nabla|v^\varepsilon|^2.
\end{equation}
Then we have for any $\varphi\in C_0^\infty(\Omega)$,
\begin{equation}\label{eq:derivative of Ju+Jv}
	\int_{\Omega}\p_t(J(u^\varepsilon)+J(v^\varepsilon))\varphi\ud\vx=\int_{\Omega}\left(\la \Hess(\varphi)\nabla u^\varepsilon,\J\nabla u^\varepsilon\ra+\la \Hess(\varphi)\nabla v^\varepsilon,\J\nabla v^\varepsilon\ra\right)\ud\vx.
\end{equation}
\subsection{Existence of vortex paths}
\begin{Lem}[Existence of vortex paths]\label{lem:existence of vortex paths}
Assume $(u_0^\varepsilon,v_0^\varepsilon)$ satisfies the conditions of Theorem \ref{thm:main}, then then there exist $T>0$ and Lipschitz paths $\ta_j,\tb_k:[0,T)\to \Omega,1\le j\le M,1\le k\le N$ such that 
\begin{equation}\label{eq:convergence of Ju(t) initial}
	J(u^\varepsilon(t))\to \frac{\pi}{1+g}\sum_{j=1}^Md_j^1\delta_{\ta_j(t)},\quad J(v^\varepsilon)\to \frac{\pi}{1+g}\sum_{j=1}^Nd_j^2\delta_{\tb_j(t)}, \quad \text{in}\  W^{-1,1}(\Omega).
\end{equation}
\end{Lem}
\begin{proof}
We denote $r_0=r( (\va,\vb) )$. From \eqref{eq:convergence of Ju0} and that $(u^\varepsilon,v^\varepsilon)$ is the solution of \eqref{eq:main}, we have that there exists $\varepsilon_0>0$ and $T^\varepsilon>0$ for $\varepsilon<\varepsilon_0$ such that for any $t<T^\varepsilon$, we have
\begin{equation}
	\left\|J(u^\varepsilon_0)-\frac{\pi}{1+g}\sum_{j=1}^Md_j^1\delta_{\va_j^0}\right\|_{W^{-1,1}(\Omega)}\le \frac{\pi}{400(1+g)}r_0, \left\|J(v^\varepsilon_0)-\frac{\pi}{1+g}\sum_{j=1}^Nd_j^2\delta_{\vb_j^0}\right\|_{W^{-1,1}(\Omega)}\le \frac{\pi}{400(1+g)}r_0,
\end{equation}
\begin{equation}\label{eq:def of Te}
	\left\|J(u^\varepsilon(t))-J(u^\varepsilon_0)\right\|_{W^{-1,1}(\Omega)}\le \frac{\pi}{400(1+g)}r_0, \left\|J(v^\varepsilon(t))-J(v^\varepsilon_0)\right\|_{W^{-1,1}(\Omega)}\le \frac{\pi}{400(1+g)}r_0.
\end{equation}
Without ambiguity, we can assume that $T^\varepsilon$ is the maximum one such that \eqref{eq:def of Te} holds for any $t<T^\varepsilon$. In particular, we have for any $1\le j\le M$, 
\begin{equation}
	\left\|J(\sqrt{1+g}u^\varepsilon)-\pi d_j^1\delta_{\va_j^0}\right\|_{W^{-1,1}(B_{r_0}(\va_j^0))}\le \frac{\pi}{200}r_0.
\end{equation}
We denote $g_-=\sqrt{(1-|g|)/(1+g)}$. If 
\begin{equation}
	\int_{B_{r_0}(\va_j^0)}e^{g_-\varepsilon}(\sqrt{1+g}u^\varepsilon)\ud\vx\le \pi\log\left(\frac{r}{\varepsilon}\right)+\gamma,
\end{equation}
then according to Theorem 1.4.3 in \cite{CollianderJerrard1999GLvorticesSE},
 we have that
\begin{equation}
	\int_{B_{r_0}(\va_j^0)}e^{g_-\varepsilon}(\sqrt{1+g}u^\varepsilon)\ud\vx\ge \pi\log\left(\frac{1}{\varepsilon}\right)-C
\end{equation}
for some $C=C(\gamma)$. Otherwise, we also have
\begin{equation}
	\int_{B_{r_0}(\va_j^0)}e^{g_-\varepsilon}(\sqrt{1+g}u^\varepsilon)\ud\vx> \pi\log\left(\frac{r}{\varepsilon}\right)+\gamma\ge \pi\log\left(\frac{1}{\varepsilon}\right)-C.
\end{equation}
Hence, we have
\begin{equation}\label{eq:weak lower bound of u}
	E^{g_-\varepsilon}(\sqrt{1+g}u^\varepsilon)\ge\pi M\log\frac{1}{\varepsilon}-C,
\end{equation}
which together with \eqref{eq:energy bound of initial data}, \eqref{eq:Egu+Egv<Euv} and the conservation law \eqref{eq:conserve energy} implies
\begin{equation}\label{eq:weak upper bound of v}
	E^{g_-\varepsilon}(\sqrt{1+g}v^\varepsilon)\le\pi N\log\frac{1}{\varepsilon}+C,
\end{equation}
Similarly,
\begin{equation}\label{eq:weak lower bound of v}
	E^{g_-\varepsilon}(\sqrt{1+g}v^\varepsilon)\ge\pi N\log\frac{1}{\varepsilon}-C,
\end{equation}
\begin{equation}\label{eq:weak upper bound of u}
	E^{g_-\varepsilon}(\sqrt{1+g}u^\varepsilon)\le\pi M\log\frac{1}{\varepsilon}+C.
\end{equation}
Hence, according to Theorem 1.4.4 in \cite{CollianderJerrard1999GLvorticesSE}, we know there exist $\va_j^\varepsilon(t)\in B_{r_0/2}(\va_j^0),\vb_j^\varepsilon(t)\in B_{r_0/2}(\vb_j^0)$ for any fixed $\varepsilon<\varepsilon_0$ and $t<T^\varepsilon$ such that 
\begin{equation}\label{eq:J is close to deltae}
	\left\|J(u^\varepsilon(t))-\frac{\pi}{1+g}\sum_{j=1}^Md_j^1\delta_{\va_j^\varepsilon(t)}\right\|_{W^{-1,1}(\Omega)}=o(1),\quad \left\|J(v^\varepsilon(t))-\frac{\pi}{1+g}\sum_{j=1}^Md_j^2\delta_{\vb_j^\varepsilon(t)}\right\|_{W^{-1,1}(\Omega)}=o(1),
\end{equation}
and for any $\rho<r_0/2$,
\begin{equation}\label{eq:weak energy bound out balls}
	\int_{\Omega_\rho(\va,\vb)}e^\varepsilon_0(u^\varepsilon,v^\varepsilon)\ud\vx\le C.
\end{equation}
Then we need to estimate $|\va_j^\varepsilon(t_2)-\va_j^\varepsilon(t_1)|$. We denote $\nu=d_j^1(\va_j^\varepsilon(t_2)-\va_j^\varepsilon(t_1))/|\va_j^\varepsilon(t_2)-\va_j^\varepsilon(t_1)|$ and find 
$\varphi\in C_0^\infty(B_{r_0}(\va_j^0))$ such that
\begin{equation}
	\varphi(\vx)=\nu\cdot(\vx-\va_j^0),\quad \text{in}\ B_{r_0/2}(\va_j^0).
\end{equation}
Then, substituting $\varphi$ to \eqref{eq:derivative of Ju+Jv}, and noting \eqref{eq:J is close to deltae} and \eqref{eq:weak energy bound out balls}, we have
\begin{align}\label{eq:ajt2-ajt1<t1-t2}
	|\va_j^\varepsilon(t_2)-\va_j^\varepsilon(t_1)|=&\int_{\Omega}d_j^1(\delta_{\va_j^\varepsilon(t_2)}-\delta_{\va_j^\varepsilon(t_1)})\varphi\ud\vx\nn 
	=&\int_{\Omega}(\sum_{k=1}^Md_k^1\delta_{\va_k^\varepsilon(t_2)}+\sum_{k=1}^Nd_k^2\delta_{\vb_k^\varepsilon(t_2)}-\sum_{k=1}^Md_k^1\delta_{\va_k^\varepsilon(t_1)}-\sum_{k=1}^Nd_k^2\delta_{\vb_k^\varepsilon(t_1)})\varphi\ud\vx\nn 
	=&\frac{1+g}{\pi}\int_{t_1}^{t_2}\int_{\Omega}\varphi\p_t(J(u^\varepsilon)+J(v^\varepsilon))\ud\vx \ud t+o(1)\nn 
	=&\frac{1+g}{\pi}\int_{t_1}^{t_2}\int_{\Omega}\left(\la \Hess(\varphi)\nabla u^\varepsilon,\J\nabla u^\varepsilon\ra+\la \Hess(\varphi)\nabla v^\varepsilon,\J\nabla v^\varepsilon\ra\right)\ud\vx\ud t+o(1)\nn
	\le &\frac{1+g}{\pi}\|\varphi\|_{C^2(\Omega\times[t_1,t_2])}\sup_{[t_1,t_2]}(\|\nabla u^\varepsilon(t)\|^2_{L^2(\Omega)}+\|\nabla v^\varepsilon(t)\|^2_{L^2(\Omega)})|t_1-t_2|+o(1)\nn
	 \le& C|t_2-t_1|+o(1).
\end{align}
Repeating the proof of Arzel\`a-Ascoli theorem with some adjustments, \eqref{eq:ajt2-ajt1<t1-t2} implies that we can find Lipschitz paths $\ta_j$ such that up to a subsequence,
\begin{equation}
	\va_j^\varepsilon(t)\to \ta_j(t).
\end{equation}
Similarly, we can also find $\tb_j$ such that 
\begin{equation}
 	\vb_j^\varepsilon(t)\to \tb_j(t).
 \end{equation} 
\end{proof}

\subsection{Convergence of the current of the solution}
Recall $d^1,d^2,\ta=(\ta_1,\cdots,\ta_M),\tb=(\tb_1,\cdots,\tb_N)$ are given in Lemma \ref{lem:existence of vortex paths}, we define $H_{d^1}(\vx;\ta(t)),H_{d^2}(\vx;\tb)$ by solving \eqref{eq:def of H initial} with $\va$ replaced by $\ta(t)$, $\tb(t)$, respectively.
\begin{equation}\label{eq:def of u* v*}
	% \begin{cases}
	u^*(\vx,t):=\frac{1}{\sqrt{1+g}}H_{d^1}(\vx;\ta(t)),\quad v^*(\vx,t):=\frac{1}{\sqrt{1+g}}H_{d^2}(\vx;\tb(t))
	% \end{cases}
\end{equation}

We then  calculate the limit of $\vj(u^\varepsilon)$ and $\vj(v^\varepsilon)$ as $\varepsilon\to 0$.

\begin{Lem}\label{lem:jtoj*}
Recall that $\ta,\tb,T$ are obtained in Lemma \ref{lem:existence of vortex paths}, and for any  $T_1\in(0,T)$, we have that as $\varepsilon\to 0$,
\begin{equation}\label{eq:vju converges to vju*}
	\vj(u^\varepsilon)\wto \vj(u^*),\quad \vj(v^\varepsilon)\wto \vj(v^*), \quad \text{in}\ L^1(\Omega\times[0,T_1]),
\end{equation}
\begin{equation}\label{eq:vju/|u| converges to vju*}
	\frac{\vj(u^\varepsilon)}{|u^\varepsilon|}\wto \sqrt{1+g}\vj(u^*),\quad \frac{\vj(v^\varepsilon)}{|v^\varepsilon|}\wto \sqrt{1+g}\vj(v^*), \quad \text{in}\ L^2_{loc}(\Omega_*(\ta(t),\tb(t))\times[0,T_1]).
\end{equation}
\end{Lem}
\begin{proof}
Combing \eqref{eq:weak upper bound of u}, \eqref{eq:convergence of Ju(t)} and Theorem 1.4.3 in \cite{CollianderJerrard1999GLvorticesSE}, we have that
\begin{equation}
	\|\vj(u^\varepsilon)\|_{L^1(\Omega)}\le C,\quad \|\vj(v^\varepsilon)\|_{L^1(\Omega)}\le C.
\end{equation}
Hence, there exist $\vj_{a},\vj_b\in L^1(\Omega\times[0,T_1])$ such that
\begin{equation}
	\vj(u^\varepsilon)\wto \vj_a,\quad \vj(v^\varepsilon)\wto \vj_b.
\end{equation}
Combining \eqref{eq:conserve mass} and the convergence of $\vj(u^\varepsilon)$, we have that for any $\varphi\in C^\infty_0(\Omega\times[0,T_1])$, 
\begin{align}\label{eq:D vja=0}
	\int_{\Omega\times[0,T_1]}\varphi\nabla\cdot \vj_a\ud \vx\ud t&=-\int_{\Omega\times[0,T_1]}\nabla \varphi\cdot \vj_a\ud\vx\ud t=-\lim_{\varepsilon\to 0} \int_{\Omega\times[0,T_1]}\nabla \varphi\cdot\vj(u^\varepsilon)\ud\vx\ud t\nn 
	&=\lim_{\varepsilon\to 0} \int_{\Omega\times[0,T_1]}\varphi\nabla \cdot\vj(u^\varepsilon)\ud\vx\ud t=\lim_{\varepsilon\to 0} \int_{\Omega\times[0,T_1]}\varphi\frac{\p |u^\varepsilon|^2}{\p t}\ud\vx\ud t\nn 
	&=-\lim_{\varepsilon\to 0} \int_{\Omega\times[0,T_1]}\frac{\p \varphi}{\p t}(|u^\varepsilon|^2-1)\ud\vx\ud t=0,
\end{align}
where the last equality holds because of \eqref{eq:weak upper bound of u}. 
Meanwhile, with a similar method of \eqref{eq:D vja=0}, \eqref{eq:convergence of Ju(t) initial} implies that
\begin{equation}\label{eq:D Jvja=delta}
	\nabla\cdot(\J\vj_a)=\frac{\pi}{1+g}\sum_{j=1}^Md_j^1\delta_{\ta_j(t)}.
\end{equation}
Combining \eqref{eq:D vja=0}, \eqref{eq:D Jvja=delta} and \eqref{eq:def of u* v*}, we have
\begin{equation}
	\nabla\cdot\left(\vj_a-\vj(u^*) \right)=0,\quad \nabla\cdot\J\left(\vj_a-\vj(u^*) \right)=0,
\end{equation}
which together with the boundary conditions $\vn\cdot\vj(u^*)=0,\vn\cdot\vj_{\va}=\lim_{\varepsilon\to 0}\vn\cdot\vj(u^\varepsilon)=0$ implies that
\begin{equation}
	\vj_a=\vj(u^*).
\end{equation}
Similarly, we can prove $\vj_{\vb}=\vj(v^*)$, which completes the proof of  \eqref{eq:vju converges to vju*}.

According to \eqref{eq:weak energy bound out balls}, we know
\begin{equation}
	\int_{\Omega_\rho(\ta,\tb)\times[0,T_1]}\left|\frac{\vj(u^\varepsilon)}{|u^\varepsilon|} \right|^2\ud\vx\ud t\le C,
\end{equation}
which implies that there exists $\vj_a^*$ such that
\begin{equation}
	\frac{\vj(u^\varepsilon)}{|u^\varepsilon|}\wto \vj_a^*.
\end{equation}
Meanwhile, we have $|u^\varepsilon|\to 1/\sqrt{1+g}$ and 
\begin{equation}
	|u^\varepsilon|\frac{\vj(u^\varepsilon)}{|u^\varepsilon|}\wto \vj(u^*),
\end{equation}
which implies that $\vj(u^\varepsilon)/|u^\varepsilon|\wto \sqrt{1+g}\vj(u^*)$. Similarly, we can obtain $\vj(v^\varepsilon)/|v^\varepsilon|\wto \sqrt{1+g}\vj(v^*)$, which gives \eqref{eq:vju/|u| converges to vju*}.
\end{proof}

We first display a lemma here, whose proof will be put off to Section \ref{sec:lower bound near vortex} since it's very lengthy. 

\begin{Lem}\label{lem:lower bound near vortex}
% Assume $J(u^\varepsilon)\wto \frac{\pi}{1+g}\delta_{\bvec{0}},J(v^\varepsilon)\wto 0$, then t
There exists $C>0,r_0>0$ such that if 
\begin{equation}
	\|J(u^\varepsilon)-\frac{\pi}{1+g}\delta_{\bvec{0}}\|_{W^{-1,1}(B_1(\bvec{0})})+\|J(v^\varepsilon)\|_{W^{-1,1}(B_1(\bvec{0}))}+\varepsilon|\log \varepsilon|<r_0,
\end{equation}
then
% \begin{equation}\label{eq:lower bound near vortex}
% 	\liminf_{\varepsilon\to 0}\left(\int_{B_1(\bvec{0})}e^\varepsilon_0(u^\varepsilon,v^\varepsilon)\ud\vx- \frac{\pi}{1+g}\log\frac{1}{\varepsilon\sqrt{1+g}}-\gamma_g\right)\ge 0.
% \end{equation}
\begin{align}\label{eq:lower bound near vortex2}
	&\int_{B_1(\bvec{0})}e^\varepsilon_0(u^\varepsilon,v^\varepsilon)\ud\vx- \frac{\pi}{1+g}\log\frac{1}{\varepsilon\sqrt{1+g}}-\gamma_g\nn 
	&\quad\ge \frac{C}{\log(\|J(u^\varepsilon)-\frac{\pi}{1+g}\delta_{\bvec{0}}\|_{W^{-1,1}(B_1(\bvec{0})})+\|J(v^\varepsilon)\|_{W^{-1,1}(B_1(\bvec{0}))}+\varepsilon|\log \varepsilon|)}.
\end{align}
\end{Lem}

\begin{Rmk}
% $J(u^\varepsilon)\wto \frac{\pi}{1+g}\delta_{\bvec{0}},J(v^\varepsilon)\wto 0$, then v
Via rescalling, Lemma \ref{lem:lower bound near vortex} gives that if
\begin{equation}
	\|J(u^\varepsilon)-\frac{\pi}{1+g}\delta_{\bvec{0}}\|_{W^{-1,1}(B_R(\bvec{0}))}/R+\|J(v^\varepsilon)\|_{W^{-1,1}(B_R(\bvec{0}))}/R<r_0,
\end{equation}
then
% With rescalling, we have 
\begin{equation}\label{eq:lower bound near vortexR}
	\int_{B_R(\bvec{0})}e^\varepsilon_0(u^\varepsilon,v^\varepsilon)\ud\vx- \frac{\pi}{1+g}\log\frac{R}{\varepsilon\sqrt{1+g}}-\gamma_g\ge C\delta_\varepsilon,
\end{equation}
with
\[
	\delta_\varepsilon=\frac{1}{\log(\|J(u^\varepsilon)-\frac{\pi}{1+g}\delta_{\bvec{0}}\|_{W^{-1,1}(B_R(\bvec{0}))}/R+\|J(v^\varepsilon)\|_{W^{-1,1}(B_R(\bvec{0}))}/R+(\varepsilon/R)|\log(R/ \varepsilon)|)}
\]
\end{Rmk}

\begin{Lem}\label{lem:e|u||v|+|j/||-j|}
We have
\begin{align}
	&\limsup_{\varepsilon\to 0}\int_{\Omega_\rho(\ta,\tb)\times[t_1,t_2]}\left(\left\|\frac{\vj(u^\varepsilon)}{|u^\varepsilon|}-\sqrt{1+g}\vj(u^*)\right\|^2+\left\|\frac{\vj(v^\varepsilon)}{|v^\varepsilon|}-\sqrt{1+g}\vj(v^*)\right\|^2+2e^\varepsilon_0(|u^\varepsilon|,|v^\varepsilon|)\right)\ud\vx\ud t\nn& \quad\le \frac{2}{1+g}\int_{t_1}^{t_2}( W_{d^1}(\va^0)+W_{d^2}(\vb^0)-W_{d_1}(\ta(t))-W_{d^2}(\tb(t)))\ud t.
\end{align}
\end{Lem}
\begin{proof}
From \eqref{eq:def of u* v*}, \eqref{eq:def of E} and (3.16) in \cite{ZhuBaoJian2023}, we have
\begin{align}
	&\int_{\Omega}e^\varepsilon_0(u^\varepsilon,v^\varepsilon)\ud\vx\nn &\quad=\sum_{j=1}^M\int_{B_\rho(\ta_j)}e^\varepsilon(u^\varepsilon,v^\varepsilon)\ud\vx+\sum_{j=1}^N\int_{B_\rho(\tb_j)}e^\varepsilon(u^\varepsilon,v^\varepsilon)\ud\vx\nn 
	&\qquad+\int_{\Omega_\rho(\ta,\tb)}\left(e^\varepsilon_0(|u^\varepsilon|,|v^\varepsilon|) +\frac{1}{2}\left|\frac{\vj(u^\varepsilon)}{|u^\varepsilon|}-\sqrt{1+g}\vj(u^*)\right|^2+\frac{1}{2}\left|\frac{\vj(v^\varepsilon)}{|v^\varepsilon|}-\sqrt{1+g}\vj(v^*)\right|^2\right)\ud\vx\nn 
	&\qquad+\int_{\Omega_\rho(\ta,\tb)}\left(\la\sqrt{1+g}\vj(u^*),\frac{\vj(u^\varepsilon)}{|u^\varepsilon|}-\sqrt{1+g}\vj(u^*)\ra+\la\sqrt{1+g}\vj(v^*),\frac{\vj(v^\varepsilon)}{|v^\varepsilon|}-\sqrt{1+g}\vj(v^*)\ra\right)\ud\vx\nn 
	&\qquad+\int_{\Omega_\rho(\ta,\tb)}\left(e^\varepsilon_0(u^*,v^*)\right)\ud\vx
\end{align}
Substituting \eqref{eq:lower bound near vortex2},\eqref{eq:conserve energy},  \eqref{eq:energy bound of initial data} and \eqref{eq:estimate of euv*} into the above equation, we have
\begin{align}
	&\int_{\Omega_\rho(\ta,\tb)}\left(e^\varepsilon_0(|u^\varepsilon|,|v^\varepsilon|) +\frac{1}{2}\left|\frac{\vj(u^\varepsilon)}{|u^\varepsilon|}-\sqrt{1+g}\vj(u^*)\right|^2+\frac{1}{2}\left|\frac{\vj(v^\varepsilon)}{|v^\varepsilon|}-\sqrt{1+g}\vj(v^*)\right|^2\right)\ud\vx\nn 
&\quad \le \frac{1}{1+g}(W_{d^1}(\va^0)+W_{d^2}(\vb^0)-W_{d_1}(\ta(t))-W_{d^2}(\tb(t)))+o(1)+O(\rho^2)\nn 
&\qquad -\int_{\Omega_\rho(\ta,\tb)}\left(\la\sqrt{1+g}\vj(u^*),\frac{\vj(u^\varepsilon)}{|u^\varepsilon|}-\sqrt{1+g}\vj(u^*)\ra+\la\sqrt{1+g}\vj(v^*),\frac{\vj(v^\varepsilon)}{|v^\varepsilon|}-\sqrt{1+g}\vj(v^*)\ra\right)\ud\vx
\end{align}
Integrating over $[t_1,t_2]$ and letting $\varepsilon\to 0$, we have
\begin{align}
	&\limsup_{\varepsilon\to 0}\int_{\Omega_\rho(\ta,\tb)\times[t_1,t_2]}\left(\left|\frac{\vj(u^\varepsilon)}{|u^\varepsilon|}-\sqrt{1+g}\vj(u^*)\right|^2+\left|\frac{\vj(v^\varepsilon)}{|v^\varepsilon|}-\sqrt{1+g}\vj(v^*)\right|^2+2e^\varepsilon_0(|u^\varepsilon|,|v^\varepsilon|)\right)\ud\vx\nn& \quad\le \frac{2}{1+g}\int_{t_1}^{t_2}( W_{d^1}(\va^0)+W_{d^2}(\vb^0)-W_{d_1}(\ta(t))-W_{d^2}(\tb(t)))\ud t+O(\rho^2).
\end{align}
Letting $\rho\to 0$ and noting that $\Omega_\rho(\ta,\tb)\subset \Omega_{\rho'}(\ta,\tb)$ for any $\rho'<\rho$, we finish the proof.
\end{proof}

\subsection{Proof of the reduced dynamical law}

In this section, we give the proof of Theorem \ref{thm:main}.

\begin{proof}
Denote the solution of \eqref{eq:MainODE} is $(\va,\vb)$ and assume $\ta,\tb$ are obtained in Lemma \ref{lem:existence of vortex paths}. Then, we only need to proof $\va=\ta,\vb=\tb$, i.e. $\zeta(t)=0$ with
\begin{equation}
	\zeta(t):=\sum_{j=1}^M|\va_j(t)-\ta_j(t)|+\sum_{j=1}^N|\vb_j(t)-\tb_j(t) |
\end{equation}
Noting \eqref{eq:MainODE}, \eqref{eq:convergence of Ju0} and \eqref{eq:convergence of Ju(t) initial}, we have $\zeta(0)=0$.

Taking derivative of $\zeta$ with respect to $t$, and noting \eqref{eq:MainODE}, we have
\begin{align}
	\dot{\zeta}(t)\le &\sum_{j=1}^M|\dot{\va}_j(t)-\dot{\ta}_j(t)|+\sum_{j=1}^N|\dot{\vb}_j(t)-\dot{\tb}_j(t) |\nn 
	\le & \sum_{j=1}^M\left|-\frac{d_j^1}{\pi}\J\nabla _{\va_j}W_{d^1}(\va)+\frac{d_j^1}{\pi}\J\nabla _{\ta_j}W_{d^1}(\ta)\right|\nn &+\sum_{j=1}^N\left|-\frac{d_j^1}{\pi}\J\nabla _{\vb_j}W_{d^1}(\vb)+\frac{d_j^1}{\pi}\J\nabla _{\tb_j}W_{d^1}(\tb)\right|\nn 
	&+\sum_{j=1}^M\left|-\frac{d_j^1}{\pi}\J\nabla _{\ta_j}W_{d^1}(\ta)-\dot{\ta}_j(t)\right|+\sum_{j=1}^N\left|-\frac{d_j^1}{\pi}\J\nabla _{\tb_j}W_{d^1}(\tb)-\dot{\tb}_j(t)\right|\nn 
	\le & C\zeta(t)+\sum_{j=1}^MA_j+\sum_{j=1}^N B_j.
\end{align}

For any $A_j$, we can find a unit-vector $\bvec{\nu}$ such that:
\begin{equation}
	A_j=d_j^1\bvec{\nu}\cdot\left(-\frac{d_j^1}{\pi}\J\nabla _{\ta_j}W_{d^1}(\ta)-\dot{\ta}_j(t)\right).
\end{equation}
Then, we can find $\varphi\in C^\infty_0(B_{r_*}(\ta_j))$ such that 
\begin{equation}\label{eq:def of varphi}
	\varphi(\vx)=\bvec{\nu}\cdot(\vx-\ta_j(t)),\quad \vx\in B_{r_*/2}(\ta_j(t)).
\end{equation}
Substituting $\varphi$ into \eqref{eq:derivative of Ju+Jv} and \eqref{eq:hessphijujv}, noting \eqref{eq:convergence of Ju(t) initial} and (4.16) in \cite{ZhuBaoJian2023}, we have
\begin{align}\label{eq:A=I+I}
A_j=&-\lim_{h\to 0}\lim_{\varepsilon\to 0}\frac{1+g}{\pi h}\int_{t}^{t+h}\int_\Omega (\la\Hess (\varphi)\nabla u^\varepsilon,\J\nabla u^\varepsilon\ra+\la\Hess (\varphi)\nabla v^\varepsilon,\J\nabla v^\varepsilon\ra)\ud\vx\ud t\nn 
&+\lim_{h\to 0}\lim_{\varepsilon\to 0}\frac{(1+g)^2}{\pi h}\int_{t}^{t+h}\int_\Omega (\la\Hess (\varphi)\vj(u^*),\J\vj(u^*)\ra+\la\Hess (\varphi)\vj(v^*),\J\vj(v^*)\ra)\ud\vx \ud t=I_1+I_2,
% = &-\lim_{h\to 0}\lim_{\varepsilon\to 0}\frac{1+g}{\pi h}\int_{t}^{t+h}\int_\Omega(\la\Hess (\varphi)\nabla |u^\varepsilon|,\J\nabla |u^\varepsilon|\ra+\la\Hess (\varphi)\nabla |v^\varepsilon|,\J\nabla |v^\varepsilon|\ra)\ud\vx\ud t\nn
% &-\lim_{h\to 0}\lim_{\varepsilon\to 0}\frac{1+g}{\pi h}\int_{t}^{t+h}\int_\Omega\left(\la\Hess (\varphi)\left(\frac{\vj(u^\varepsilon)}{|u^\varepsilon|}-\sqrt{1+g}\vj(u^*)\right),\J\left(\frac{\vj(u^\varepsilon)}{|u^\varepsilon|}-\sqrt{1+g}\vj(u^*)\right)\ra \right.\nn 
% &\quad\quad +\left.\la\Hess (\varphi)\left(\frac{\vj(v^\varepsilon)}{|v^\varepsilon|}-\sqrt{1+g}\vj(v^*)\right),\J\left(\frac{\vj(v^\varepsilon)}{|v^\varepsilon|}-\sqrt{1+g}\vj(v^*)\right)\ra \right)\ud\vx \ud t\nn
% &-\lim_{h\to 0}\lim_{\varepsilon\to 0}\frac{1+g}{\pi h}\int_{t}^{t+h}\int_\Omega\left(\la\Hess (\varphi)\left(\frac{\vj(u^\varepsilon)}{|u^\varepsilon|}-\sqrt{1+g}\vj(u^*)\right),\sqrt{1+g}\J\vj(u^*)\ra \right.\nn 
% &\quad\quad +\left.\la\Hess (\varphi)\left(\frac{\vj(v^\varepsilon)}{|v^\varepsilon|}-\sqrt{1+g}\vj(v^*)\right),\sqrt{1+g}\J\vj(v^*)\ra \right)\ud\vx \ud t\nn
% &-\lim_{h\to 0}\lim_{\varepsilon\to 0}\frac{1+g}{\pi h}\int_{t}^{t+h}\int_\Omega\left(\la\sqrt{1+g}\Hess (\varphi)\vj(u^*),\J\left(\frac{\vj(u^\varepsilon)}{|u^\varepsilon|}-\sqrt{1+g}\vj(u^*)\right)\ra \right.\nn 
% &\quad\quad +\left.\la\sqrt{1+g}\Hess (\varphi)\vj(v^*),\J\left(\frac{\vj(v^\varepsilon)}{|v^\varepsilon|}-\sqrt{1+g}\vj(v^*)\right)\ra \right)\ud\vx \ud t\nn
% \le &C( W_{d^1}(\va^0)+W_{d^2}(\vb^0)-W_{d_1}(\ta(t))-W_{d^2}(\tb(t)))\le C\zeta(t).
\end{align}
where
\begin{align}
	I_1= &-\lim_{h\to 0}\lim_{\varepsilon\to 0}\frac{1+g}{\pi h}\int_{t}^{t+h}\int_\Omega(\la\Hess (\varphi)\nabla |u^\varepsilon|,\J\nabla |u^\varepsilon|\ra+\la\Hess (\varphi)\nabla |v^\varepsilon|,\J\nabla |v^\varepsilon|\ra)\ud\vx\ud t\nn
&-\lim_{h\to 0}\lim_{\varepsilon\to 0}\frac{1+g}{\pi h}\int_{t}^{t+h}\int_\Omega\left(\la\Hess (\varphi)\left(\frac{\vj(u^\varepsilon)}{|u^\varepsilon|}-\sqrt{1+g}\vj(u^*)\right),\J\left(\frac{\vj(u^\varepsilon)}{|u^\varepsilon|}-\sqrt{1+g}\vj(u^*)\right)\ra \right.\nn 
&\quad\quad +\left.\la\Hess (\varphi)\left(\frac{\vj(v^\varepsilon)}{|v^\varepsilon|}-\sqrt{1+g}\vj(v^*)\right),\J\left(\frac{\vj(v^\varepsilon)}{|v^\varepsilon|}-\sqrt{1+g}\vj(v^*)\right)\ra \right)\ud\vx \ud t,
\end{align}
\begin{align}
I_2=&-\lim_{h\to 0}\lim_{\varepsilon\to 0}\frac{1+g}{\pi h}\int_{t}^{t+h}\int_\Omega\left(\la\Hess (\varphi)\left(\frac{\vj(u^\varepsilon)}{|u^\varepsilon|}-\sqrt{1+g}\vj(u^*)\right),\sqrt{1+g}\J\vj(u^*)\ra \right.\nn 
&\quad\quad +\left.\la\Hess (\varphi)\left(\frac{\vj(v^\varepsilon)}{|v^\varepsilon|}-\sqrt{1+g}\vj(v^*)\right),\sqrt{1+g}\J\vj(v^*)\ra \right)\ud\vx \ud t\nn
&-\lim_{h\to 0}\lim_{\varepsilon\to 0}\frac{1+g}{\pi h}\int_{t}^{t+h}\int_\Omega\left(\la\sqrt{1+g}\Hess (\varphi)\vj(u^*),\J\left(\frac{\vj(u^\varepsilon)}{|u^\varepsilon|}-\sqrt{1+g}\vj(u^*)\right)\ra \right.\nn 
&\quad\quad +\left.\la\sqrt{1+g}\Hess (\varphi)\vj(v^*),\J\left(\frac{\vj(v^\varepsilon)}{|v^\varepsilon|}-\sqrt{1+g}\vj(v^*)\right)\ra \right)\ud\vx \ud t
\end{align}
Lemma \ref{lem:e|u||v|+|j/||-j|} implies that 
\begin{equation}\label{eq:est of I1 1}
	|I_1|\le C( W_{d^1}(\va^0)+W_{d^2}(\vb^0)-W_{d_1}(\ta(t))-W_{d^2}(\tb(t))).
\end{equation}
Taking derivative of $W_{d^1}(\va)+W_{d^2}(\vb)$ with respect to $t$ and noting \eqref{eq:MainODE}, we will obtain that $W_{d^1}(\va(t))+W_{d^2}(\vb(t))=W_{d^1}(\va^0)+W_{d^2}(\vb^0)$. Hence \eqref{eq:est of I1 1} gives
\begin{equation}\label{eq:est of I1 2}
	|I_1|\le C( W_{d^1}(\va(t))+W_{d^2}(\vb(t))-W_{d_1}(\ta(t))-W_{d^2}(\tb(t)))\le C\zeta(t).
\end{equation}
Noting $\Hess \varphi(\vx)=0 $ in $\Omega_{r_*/2}(\ta(t,\tb(t)))$, Lemma \ref{lem:jtoj*} implies $I_2=0$. Hence \eqref{eq:A=I+I} and \eqref{eq:est of I1 2} imply 
\begin{equation}
	A_j\le C\zeta(t).
\end{equation}
Similarly, we have $B_j\le C\zeta(t)$ and hence $\dot{\zeta}(t)\le C\zeta(t)$. Then we have $\zeta(t)\equiv 0$, which finishes the proof.
\end{proof}
\section{The lower bound of energy in a neighborhood of a vortex}\label{sec:lower bound near vortex}

We give a lemma on the winding number of $u,v$: 
\begin{Lem}\label{lem:number of r such that deg is 1}
Assume \begin{equation}
	\left\|J(u)-\frac{\pi}{1+g}\delta_{\bvec{0}} \right\| +\|J(v)\|\le \frac{1}{4},
\end{equation}
\begin{equation}\label{eq:euv B1}
	\int_{B_1(\bvec{0})}e^\varepsilon_0(u,v)\ud\vx\le 2\frac{\pi}{1+g}\log\frac{1}{\varepsilon}.
\end{equation}
Then
\[
	\mathcal{L}^1\left(\left\{r\in[0,1]|deg(u,\p B_r(\bvec{0}))=1,deg(v,\p B_r(\bvec{0}))=0,|u|,|v|>\frac{1}{2}\right\}\right)\ge 1-C\left(\left\|J(u)-\frac{\pi}{1+g}\delta_{\bvec{0}} \right\| +\|J(v)\|+\varepsilon|\log\varepsilon|\right).
\]
\end{Lem}
\begin{proof}
Define
\begin{equation}
	J'(u)=\zeta(|u|)J(u),\quad J'v=\zeta(|v|)J(v),
\end{equation}
where $\zeta:[0,\infty)\to[0,\infty)$ satisfies $supp \zeta=[0,1/2]$ and $\int_{\R^2}\zeta(|\vx|)\ud\vx =\pi(1+g)$.

Then, similar to Lemma 2.1 and (2.17) in \cite{Jerrard2007RefinedJacobian}, we have
\begin{equation}\label{eq:J'u-Ju}
	\|J'(u)-J(u)\|\le C\|\nabla u\|_{L^2(B_1)}\left\|\frac{1}{1+g}-|u|^2\right\|_{L^2(B_1)},\quad \|J'(v)-J(v)\|\le C\|\nabla v\|_{L^2(B_1)}\left\|\frac{1}{1+g}-|v|^2\right\|_{L^2(B_1)},
\end{equation}
\begin{equation}\label{eq:degupV=intVJ'u}
	deg(u;\p V) =\frac{1+g}{\pi}\int_V J'(u)\ud\vx, deg(v;\p V) =\frac{1+g}{\pi}\int_V J'(v)\ud\vx, \ \text{if}\ |u|,|v|\ge \frac{1}{2}\ \text{on}\ \p V.
\end{equation}
Noting \eqref{eq:Egu+Egv<Euv} and \eqref{eq:euv B1}, we have 
\begin{equation}\label{eq:Egu Egv bounded}
	\int_{B_1}e^{g_-\varepsilon}(u)\ud\vx+\int_{B_1}e^{g_-\varepsilon}(v)\ud\vx \le 2\pi|\log\varepsilon|,
\end{equation}
which together with Lemma 3.6 in \cite{Jerrard2007RefinedJacobian} implies that 
\begin{equation}\label{eq:S small }
	\mathcal{L}^1(S)\le C\varepsilon|\log \varepsilon|,
\end{equation}
with
\begin{equation}
	S=\{r\in[0,1]:|u|\le 1/2\ \text{or}\ |v|\le 1/2 \ \text{on}\ B_r\}.
\end{equation}
Define 
\begin{equation}\label{eq:dj=intJ'}
	d_1(r)=\frac{1+g}{\pi}\int_{B_r}J'(u)\ud\vx,\quad d_2(r)=\frac{1+g}{\pi}\int_{B_r}J'(v)\ud\vx,
\end{equation}
then it follows from \eqref{eq:degupV=intVJ'u} that
\begin{equation}\label{eq:dj in Z}
	d_1(r)=deg(u;\p B_r)\in \Z, d_2(r)=deg(v;\p B_r)\in\Z,\quad r\notin S.
\end{equation}
Meanwhile, \eqref{eq:J'u-Ju} and \eqref{eq:Egu Egv bounded} imply that
\begin{align}\label{eq:J'u-delta0}
	\left\|\frac{\pi}{1+g}\delta_{\bvec{0}}-J'(u)\right\|\le &\left\|\frac{\pi}{1+g}\delta_{\bvec{0}}-J(u)\right\|+\|J'(u)-J(u)\|	\le \left\|\frac{\pi}{1+g}\delta_{\bvec{0}}-J(u)\right\|+C\varepsilon|\log\varepsilon|,
\end{align}
\begin{equation}
	\|J'(v)\|\le \|J(v)\|+C\varepsilon|\log\varepsilon|,
\end{equation}
Consider $\varphi(\vx)=f(|\vx|)$ with $f:\R\to\R$ satisfying $f(1)=0,|f'|\le 1$. Then \eqref{eq:dj=intJ'} gives
\begin{equation}
	\int_{B_1}\varphi J'(u)\ud\vx=\int_0^1f(r)\int_{\p B_r}J'(u)\ud \sigma d r=\frac{\pi}{1+g}\int_0^1f(r)d'_1(r) d r=-\frac{\pi}{1+g}\int_0^1f'(r)d(r)\ud r,
\end{equation}
which together with \eqref{eq:J'u-delta0} gives
\begin{equation}\label{eq:f' 1-d1}
	\frac{\pi}{1+g}\int_0^1 f'(r)(1-d_1(r))\ud r\le \left\|\frac{\pi}{1+g}\delta_{\bvec{0}}-J(u)\right\|+C\varepsilon|\log\varepsilon|.
\end{equation}
Since \eqref{eq:f' 1-d1} is valid for any $f$ such that $f(1)=0, |f'|\le 1$, we have 
\begin{equation}
	\frac{\pi}{1+g}\int_0^1|1-d_1(r)|\ud r\le \left\|\frac{\pi}{1+g}\delta_{\bvec{0}}-J(u)\right\|+C\varepsilon|\log\varepsilon|,
\end{equation}
which together with \eqref{eq:dj in Z} and \eqref{eq:S small } implies that
\begin{equation}\label{eq:degu ne 1}
	\mathcal{L}^1\left(\{r\in [0,1]| deg(u,\p B_r(\bvec{0}))\ne 1\}\right)\le C\left(\left\|\frac{\pi}{1+g}\delta_{\bvec{0}}-J(u)\right\|+\varepsilon|\log\varepsilon|\right).
\end{equation}
Similarly,
\begin{equation}\label{eq:degv ne 0}
	\mathcal{L}^1\left(\{r\in [0,1]| deg(v,\p B_r(\bvec{0}))\ne 0\}\right)\le C\left(\left\|J(v)\right\|+\varepsilon|\log\varepsilon|\right).
\end{equation}
We then complete the proof via combining \eqref{eq:S small }, \eqref{eq:degu ne 1} and \eqref{eq:degv ne 0}.
\end{proof}

Here, we give the proof of Lemma \ref{lem:lower bound near vortex}.

\begin{proof}
Let $\delta_\varepsilon=\|J(u^\varepsilon)-\frac{\pi}{1+g}\delta_{\bvec{0}}\|_{W^{-1,1}(B_1(\bvec{0}))})+\|J(v^\varepsilon)\|_{W^{-1,1}(B_1(\bvec{0}))}+\varepsilon|\log \varepsilon|$. Without loss of generality, we can assume 
\[
	\int_{B_1(\bvec{0})}e^\varepsilon_0(u^\varepsilon,v^\varepsilon)\ud\vx\le \frac{\pi}{1+g}\log\frac{1}{\varepsilon\sqrt{1+g}}+\gamma_g.
\]
Otherwise, we can take $C=0$.

On the other hand, the fact that $\|J(u^\varepsilon)-\frac{\pi}{1+g}\delta_{\bvec{0}}\|_{W^{-1,1}(B_1(\bvec{0}))}$ implies that for any $\varepsilon<<r_0<1$,
\[
	\int_{B_{r_0}(\bvec{0})}e^\varepsilon_0(u^\varepsilon,v^\varepsilon)\ud\vx\ge \frac{\pi}{1+g}\int_{B_1(\bvec{0})}e^{g_-\varepsilon}(\sqrt{1+g}u^\varepsilon)\ud\vx\ge \frac{\pi}{1+g}\log\frac{r_0}{\varepsilon}-C.
\]
Hence, there exists $C_1$ such that for any $\varepsilon<<r_0<1$, we have
\begin{equation}
	\int_{B_{1/2}\setminus B_{r_0}(\bvec{0})}e^\varepsilon_0(u^\varepsilon,v^\varepsilon)\ud\vx\le \frac{\pi}{1+g}\log\frac{1/2}{r_0}+C_1.
\end{equation}
Then we take $r_0=\delta_\varepsilon^{1/4}/2$, which implies that 
\begin{equation}\label{eq:upper bound of energy in B-Br}
	\int_{B_{1/2}\setminus B_{r_0}(\bvec{0})}e^\varepsilon_0(u^\varepsilon,v^\varepsilon)\ud\vx\le \frac{\pi}{1+g}\log\frac{1/2}{r_0}+C_1\le\frac{\pi+\tilde{\delta}_\varepsilon}{1+g}\log\frac{1/2}{r_0},
\end{equation}
with $\tilde{\delta}_\varepsilon=-4(1+g)C_1/\log\delta_\varepsilon$.

Let 
\[
	T_1=\left\{r\in[r_0,1/2]|\int_{\p B_r(\bvec{0})}e^\varepsilon_0(u^\varepsilon,v^\varepsilon)\le \frac{\pi+2\tilde{\delta}_\varepsilon}{(1+g)r} \right\},
\]
and
\[
	T_2=\left\{r\in[r_0,1/2]| deg(u^\varepsilon,\p B_r(\bvec{0}))=1,deg(v^\varepsilon,\p B_r(\bvec{0}))=0,|u^\varepsilon|,|v^\varepsilon|>\frac{1}{2}\right\},
\]
Noting that 
\begin{align}\label{eq:lower bound of T1}
\int_{B_{1/2}\setminus B_{r_0}(\bvec{0})}e^\varepsilon_0(u^\varepsilon,v^\varepsilon)d\vx &=\int_{r_0}^{1/2}\int_{\p B_r(\bvec{0})}e^\varepsilon_0(u^\varepsilon,v^\varepsilon)\ud\sigma\ud r\nn 
&\ge \int_{[r_0,1/2]\setminus T_1}\int_{\p B_r(\bvec{0})}e^\varepsilon_0(u^\varepsilon,v^\varepsilon)\ud\sigma\ud r\nn 
&\ge \int_{[r_0,1/2]\setminus T_1}\frac{\pi+2\tilde{\delta}_\varepsilon}{(1+g)r}\ud r\nn 
&\ge \int_{r_0+\mathcal{L}^1(T_1)}^{1/2}\frac{\pi+2\tilde{\delta}_\varepsilon}{(1+g)r}\ud r\ge \frac{\pi+2\tilde{\delta}_\varepsilon}{1+g}\log\frac{1/2}{r_0+\mathcal{L}^1(T_1)}.
\end{align}
Comparing \eqref{eq:lower bound of T1} and \eqref{eq:upper bound of energy in B-Br}, we have 
\begin{equation}
	\mathcal{L}^1(T_1)\ge \left(\frac{1}{2}\right)^\frac{\tilde{\delta}_\varepsilon}{\pi+2\tilde{\delta}_\varepsilon}r_0^\frac{\pi+\tilde{\delta}_\varepsilon}{\pi+2\tilde{\delta}_\varepsilon}-r_0.
\end{equation}
Meanwhile, Lemma \ref{lem:number of r such that deg is 1} implies 
\begin{equation}
	\mathcal{L}^1(T_2)\ge\frac{1}{2}-r_0-C\delta_\varepsilon.
\end{equation}
Then we have that for $\delta_\varepsilon$ small enough,
\begin{align*}
	\mathcal{L}^1(T_1\cap T_2)&\ge \mathcal{L}^1(T_1)+ \mathcal{L}^1(T_2)-\frac{1}{2}+r_0\ge\\
	&\ge \left(\frac{1}{2}\right)^\frac{\tilde{\delta}_\varepsilon}{\pi+2\tilde{\delta}_\varepsilon}r_0^\frac{\pi+\tilde{\delta}_\varepsilon}{\pi+2\tilde{\delta}_\varepsilon}-r_0 -C\delta_\varepsilon=r_0\left(e^{\frac{\tilde{\delta}_\varepsilon}{\pi+\tilde{\delta}_\varepsilon}\log\frac{1}{2r_0}}-1\right)-C\delta_\varepsilon\\
	&= \frac{\delta_\varepsilon^{1/4}}{2}(e^\frac{(1+g)C_1}{\pi+2\tilde{\delta}_\varepsilon}-1)-C\delta_\varepsilon>0.
\end{align*}

Then we can find $r_\varepsilon\in T_1\cap T_2$, i.e.
\begin{equation}
	\int_{\p B_{r_\varepsilon}(\bvec{0})}e^\varepsilon_0(u^\varepsilon,v^\varepsilon)\le \frac{\pi+2\tilde{\delta}_\varepsilon}{(1+g)r_\varepsilon}, \quad deg(u^\varepsilon,\p B_{r_\varepsilon}(\bvec{0}))=1,deg(v^\varepsilon,\p B_{r_\varepsilon}(\bvec{0}))=0,|u^\varepsilon|,|v^\varepsilon|>\frac{1}{2}.
\end{equation}
Then \eqref{eq:Egu+Egv<Euv} gives 
\begin{equation}\label{eq:lose bound of energy on boundary}
	\int_{\p B_{r_\varepsilon}(\bvec{0})}e^{g_-\varepsilon}(\sqrt{1+g}u^\varepsilon)\ud\sigma+\int_{\p B_{r_\varepsilon}(\bvec{0})}e^{g_-\varepsilon}(\sqrt{1+g}v^\varepsilon)\ud\sigma\le \frac{\pi+2\tilde{\delta}_\varepsilon}{r_\varepsilon}
\end{equation}

From Lemma 3.1.3  in \cite{CollianderJerrard1999GLvorticesSE}, we see that 
\begin{equation}
	\|1-\sqrt{1+g}|u^\varepsilon|\|_{L^\infty(\p B_{r_\varepsilon}(\bvec{0}))}\le (C\varepsilon/r_\varepsilon)^{1/2}\le C\varepsilon^{3/8},\quad \|1-\sqrt{1+g}|v^\varepsilon|\|_{L^\infty(\p B_{r_\varepsilon}(\bvec{0}))}\le (C\varepsilon/r_\varepsilon)^{1/2}\le C\varepsilon^{3/8}.
\end{equation}

Let $\tilde{r}_\varepsilon=r_\varepsilon+\varepsilon^{1/4}$. We define $\tilde{u}^\varepsilon,\tilde{v}^\varepsilon$ on $B_{\tilde{r}_\varepsilon}(\bvec{0})$ via
\begin{equation}
	\tilde{u}^\varepsilon=\left\{\begin{array}{ll}u^\varepsilon,&\vx\in B_{r_\varepsilon}(\bvec{0}),\\\frac{u^\varepsilon(r_\varepsilon\vx/|\vx|)}{|u^\varepsilon(r_\varepsilon\vx/|\vx|)|}\sqrt{|u^\varepsilon(r_\varepsilon\vx/|\vx|)|^2\frac{\tilde{r}_\varepsilon-|\vx|}{\tilde{r}_\varepsilon-r_\varepsilon}+ \frac{1}{1+g}\frac{|\vx|-r_\varepsilon}{\tilde{r}_\varepsilon-r_\varepsilon}},&\vx\in B_{\tilde{r}_\varepsilon}\setminus B_{r_\varepsilon}(\bvec{0}).
	 \end{array} \right.
\end{equation}
\begin{equation}
	\tilde{v}^\varepsilon\left\{\begin{array}{ll}v^\varepsilon,&\vx\in B_{r_\varepsilon}(\bvec{0}),\\\frac{v^\varepsilon(r_\varepsilon\vx/|\vx|)}{|v^\varepsilon(r_\varepsilon\vx/|\vx|)|}\sqrt{|v^\varepsilon(r_\varepsilon\vx/|\vx|)|^2\frac{\tilde{r}_\varepsilon-|\vx|}{\tilde{r}_\varepsilon-r_\varepsilon}+\frac{1}{1+g} \frac{|\vx|-r_\varepsilon}{\tilde{r}_\varepsilon-r_\varepsilon}},&\vx\in B_{\tilde{r}_\varepsilon}\setminus B_{r_\varepsilon}(\bvec{0}).
	 \end{array} \right.
\end{equation}

Then we can check that

\begin{equation}\label{eq:Etue Btre-Bre}
	\int_{B_{\tilde{r}_\varepsilon}\setminus B_{r_\varepsilon}(\bvec{0})}e^\varepsilon_0(\tilde{u}^\varepsilon,\tilde{v}^\varepsilon)\le \frac{\pi}{1+g}\log\frac{\tilde{r}_\varepsilon}{r_\varepsilon}+\frac{4\log 2}{1+g}\tilde{\delta}_\varepsilon ,
\end{equation}
\begin{equation}\label{eq:upper bound of Etu pBtre}
	\int_{\p B_{\tilde{r}_\varepsilon}(\bvec{0})}e^\varepsilon_0(\tilde{u}^\varepsilon,\tilde{v}^\varepsilon)\ud \sigma\le \frac{\pi+4\tilde{\delta}_\varepsilon}{(1+g)\tilde{r}_\varepsilon}.
\end{equation}

Since $deg(\tilde{u}^\varepsilon,\p B_{\tilde{r}_\varepsilon}(\bvec{0}))=deg(\tilde{u}^\varepsilon,\p B_{{r}_\varepsilon}(\bvec{0}))=1,deg(\tilde{v}^\varepsilon,\p B_{\tilde{r}_\varepsilon}(\bvec{0}))=deg(\tilde{v}^\varepsilon,\p B_{{r}_\varepsilon}(\bvec{0}))=0$, we can find $\alpha_1,\alpha_2\in\R,\varphi_1,\varphi_2\in H^1(\p B_{\tilde{r}_\varepsilon}(\bvec{0}))$ such that
\begin{equation}
	\tilde{u}^\varepsilon(\vx)=\frac{1}{\sqrt{1+g}}e^{i(\alpha_1+\theta(\vx)+\varphi_1(\vx))},\quad \tilde{v}^\varepsilon(\vx)=\frac{1}{\sqrt{1+g}}e^{i(\alpha_2+\varphi_2(\vx))},\quad \vx\in\p B_{\tilde{r}_\varepsilon}(\bvec{0}),
\end{equation}
with
\begin{equation}
	e^{i\theta(\vx)}=\frac{x+iy}{|\vx|},\quad \int_{\p B_{\tilde{r}_\varepsilon}(\bvec{0})}\varphi_1(\vx)\ud \sigma=0=\int_{\p B_{\tilde{r}_\varepsilon}(\bvec{0})}\varphi_2(\vx)\ud \sigma.
\end{equation}
Then we can extend $\tilde{u}^\varepsilon,\tilde{v}^\varepsilon$ to $B_{1}(\bvec{0})$ by defining 
\begin{equation}
 	\tilde{u}^\varepsilon(\vx)=\frac{1}{\sqrt{1+g}}e^{i(\alpha_1+\theta(\vx)+\lambda(|\vx|)\varphi_1(\vx))},\quad \tilde{v}^\varepsilon(\vx)=\frac{1}{\sqrt{1+g}}e^{i(\alpha_2+\lambda(|\vx|)\varphi_2(\vx))}, \quad \vx\in B_1\setminus B_{\tilde{r}_\varepsilon}(\bvec{0}),
 \end{equation} 
 where 
 \begin{equation}
 	\varphi_1(\vx)=\varphi_1(\tilde{r}_\varepsilon\vx/|\vx|),\quad \varphi_2(\vx)=\varphi_2(\tilde{r}_\varepsilon\vx/|\vx|),\quad \vx\in B_1\setminus B_{\tilde{r}_\varepsilon}(\bvec{0}),
 \end{equation}
 and 
 \begin{equation}
 	\lambda(r)=\left\{\begin{array}{ll}-\frac{1}{\tilde{r}_\varepsilon}(r-2\tilde{r}_\varepsilon),&r\in[\tilde{r}_\varepsilon,2\tilde{r}_\varepsilon],\\0,&r\in[2\tilde{r}_\varepsilon,1]\end{array}\right.
 \end{equation}
We have
\begin{align}\label{eq:energy of ut B1-Btre}
	\int_{B_1\setminus B_{\tilde{r}_\varepsilon}(\bvec{0})}e^\varepsilon_0(\tilde{u}^\varepsilon,\tilde{v}^\varepsilon)\ud\vx=&\frac{1}{2(1+g)}\int_{B_1\setminus B_{\tilde{r}_\varepsilon}(\bvec{0})}\left(\left|\nabla \theta+\lambda'(|\vx|)\varphi_1(\vx)\frac{\vx}{|\vx|}+\lambda(|\vx|)\nabla \varphi_1(\vx)\right|^2\right.\nn&\qquad\qquad\qquad\qquad+\left.\left|\lambda'(|\vx|)\varphi_2(\vx)\frac{\vx}{|\vx|}+\lambda(|\vx|)\nabla \varphi_2(\vx)\right|^2 \right)\ud\vx
\end{align}
Noting that $\nabla \theta(\vx)\cdot\vx=0,\nabla \varphi_1(\vx)\cdot\vx=0,\nabla \varphi_2(\vx)\cdot \vx=0$, and for any $s\in[\tilde{r}_\varepsilon,1]$,
\begin{equation}
	\int_{\p B_s(\bvec{0})}\nabla \theta(\vx)\cdot\nabla \varphi_1(\vx)\ud\sigma=\int_{\p B_s(\bvec{0})}\frac{\p}{\p \vt}\theta(\vx)\cdot\frac{\p}{\p \vt}\varphi_1(\vx)\ud\sigma=\frac{1}{s}\int_{\p B_s(\bvec{0})}\frac{\p}{\p \vt}\varphi_1(\vx)\ud\sigma=0,
\end{equation}
\eqref{eq:energy of ut B1-Btre} gives
\begin{align}\label{eq:energy ut B1-Btre 2}
	\int_{B_1\setminus B_{\tilde{r}_\varepsilon}(\bvec{0})}e^\varepsilon_0(\tilde{u}^\varepsilon,\tilde{v}^\varepsilon)\ud\vx=&\frac{1}{2(1+g)}\int_{B_1\setminus B_{\tilde{r}_\varepsilon}(\bvec{0})}\left(\frac{1}{|\vx|^2}+|\lambda'(|\vx|)|^2\varphi_1^2(\vx)+|\lambda(|\vx|)|^2\left|\frac{\p}{\p \vt}\varphi_1(\vx)\right|^2\right.\nn &\qquad\qquad\qquad\qquad\left.+ |\lambda'(|\vx|)|^2\varphi_2^2(\vx)+|\lambda(|\vx|)|^2\left|\frac{\p}{\p \vt}\varphi_2(\vx)\right|^2\right)\ud\vx\nn 
	=&\frac{1}{2(1+g)}\int_{\tilde{r}_\varepsilon}^1\left(\frac{2\pi}{r}+|\lambda'(r)|^2\int_{\p B_r(\bvec{0})}(\varphi_1^2(\vx)+\varphi_2^2(\vx))\ud\sigma\right.\nn 
	&\qquad\qquad\qquad\qquad+\left.|\lambda(|\vx|)|^2\int_{\p B_r(\bvec{0})}\left(\left|\frac{\p}{\p \vt}\varphi_1(\vx)\right|^2+\left|\frac{\p}{\p \vt}\varphi_2(\vx)\right|^2\right)\ud\sigma\right.\ud r.
\end{align}
Meanwhile, 
\begin{equation}
	\int_{\p B_{\tilde{r}_\varepsilon}(\bvec{0})}e^\varepsilon_0(\tilde{u}^\varepsilon,\tilde{v}^\varepsilon)\ud\sigma\ge \frac{1}{1+g}\left(\frac{\pi}{\tilde{r}_\varepsilon}+\int_{\p B_{\tilde{r}_\varepsilon}(\bvec{0})}\left(\left|\frac{\p}{\p \vt}\varphi_1(\vx)\right|^2+\left|\frac{\p}{\p \vt}\varphi_2(\vx)\right|^2\right)\ud\sigma\right),
\end{equation}
which together with \eqref{eq:upper bound of Etu pBtre} implies
\begin{equation}\label{eq:df2 pBtre}
	\int_{\p B_{\tilde{r}_\varepsilon}(\bvec{0})}\left(\left|\frac{\p}{\p \vt}\varphi_1(\vx)\right|^2+\left|\frac{\p}{\p \vt}\varphi_2(\vx)\right|^2\right)\ud\sigma \le \frac{4\tilde{\delta}_\varepsilon}{\tilde{r}_\varepsilon}.
\end{equation}
With Pioncare's inequality, we have
\begin{equation}\label{eq:f2 pBtre}
	\int_{\p B_{\tilde{r}_\varepsilon}(\bvec{0})}\left(\left|\varphi_1(\vx)\right|^2+\left|\varphi_2(\vx)\right|^2\right)\ud\sigma \le 4\tilde{r}_\varepsilon\tilde{\delta}_\varepsilon.
\end{equation}
Noting the homogeneity of $\varphi_1,\varphi_2$, \eqref{eq:df2 pBtre} and \eqref{eq:f2 pBtre} imply
\begin{equation}\label{eq:f2 pBr}
	\int_{\p B_{r}(\bvec{0})}\left(\left|\frac{\p}{\p \vt}\varphi_1(\vx)\right|^2+\left|\frac{\p}{\p \vt}\varphi_2(\vx)\right|^2\right)\ud\sigma \le \frac{4\tilde{\delta}_\varepsilon}{r},\quad \int_{\p B_{r}(\bvec{0})}\left(\left|\varphi_1(\vx)\right|^2+\left|\varphi_2(\vx)\right|^2\right)\ud\sigma \le 4r{\delta}_\varepsilon,\quad r\in [\tilde{r}_\varepsilon,1].
\end{equation}
Substituting \eqref{eq:f2 pBr} to \eqref{eq:energy ut B1-Btre 2}, and noting the definition of $\lambda$ we have
\begin{equation}\label{eq:energy tue B1-Btre 3}
	\int_{B_1\setminus B_{\tilde{r}_\varepsilon}(\bvec{0})}e^\varepsilon_0(\tilde{u}^\varepsilon,\tilde{v}^\varepsilon)\ud\vx\le \frac{\pi}{(1+g)}\log\frac{1}{\tilde{r}_\varepsilon}+\frac{2\tilde{\delta}_\varepsilon}{1+g}\int_{\tilde{r}_\varepsilon}^{2\tilde{r}_\varepsilon}\left(\frac{1}{r}+\frac{1}{\tilde{r}_\varepsilon^2}r \right)d r \le \frac{\pi}{(1+g)}\log\frac{1}{\tilde{r}_\varepsilon}+\frac{\tilde{\delta}_\varepsilon}{1+g}(2\log 2+3).
\end{equation}
Combining \eqref{eq:energy tue B1-Btre 3} and \eqref{eq:Etue Btre-Bre}, we have
\begin{equation}\label{eq:Etue B1}
	\int_{B_1(\bvec{0})}e^\varepsilon_0(\tilde{u}^\varepsilon,\tilde{v}^\varepsilon)\ud \vx\le \int_{B_{r_\varepsilon(\bvec{0})}}e^\varepsilon_0(u^\varepsilon,{v}^\varepsilon)\ud \vx+\frac{\pi}{(1+g)}\log\frac{1}{\tilde{r}_\varepsilon}+C\tilde{\delta}_\varepsilon.
\end{equation}

On the other hand, 
\begin{equation}
	(1+g)e_0^\varepsilon(u^\varepsilon,v^\varepsilon)\ge e^{g_-\varepsilon}(\sqrt{1+g}u^\varepsilon)+e^{g_-\varepsilon}(\sqrt{1+g}v^\varepsilon)\ge e^{g_-\varepsilon}(\sqrt{1+g}u^\varepsilon),
\end{equation}
hence
\begin{equation}
	(1+g)\int_{B_1\setminus B_{r_\varepsilon}(\bvec{0})}e_0^\varepsilon(u^\varepsilon,v^\varepsilon)\ud\vx\ge\int_{B_1\setminus B_{r_\varepsilon}(\bvec{0})} e^{g_-\varepsilon}(\sqrt{1+g}u^\varepsilon)\ud\vx.
\end{equation}
Let 
\[
	T_\varepsilon=\left\{r\in[r_\varepsilon,1]| deg(u^\varepsilon,\p B_r(\bvec{0}))=1,deg(v^\varepsilon,\p B_r(\bvec{0}))=0,|u^\varepsilon|,|v^\varepsilon|>\frac{1}{2}\right\}.
\]
Then we have $\mathcal{L}^1(T_\varepsilon)\ge 1-r_\varepsilon-C\delta_\varepsilon$ and 
\begin{align}\label{eq:Eue B1-Bre}
(1+g)\int_{B_1\setminus B_{r_\varepsilon}(\bvec{0})}e_0^\varepsilon(u^\varepsilon,v^\varepsilon)\ud\vx\ge&\int_{B_1\setminus B_{r_\varepsilon}(\bvec{0})} e^{g_-\varepsilon}(\sqrt{1+g}u^\varepsilon)\ud\vx\nn 
\ge &\int_{r\in T_\varepsilon}\int_{\p B_r(\bvec{0})}e^{g_-\varepsilon}(\sqrt{1+g}u^\varepsilon)\ud\sigma\ud r\nn
\ge&\int_{r\in T_\varepsilon}\min\{\lambda^\varepsilon(r),c_0/\varepsilon\}\ud r,
\end{align}
where we used Lemma 3.1.4 in \cite{CollianderJerrard1999GLvorticesSE} and 
\begin{equation}
	\lambda^\varepsilon(r)=\min_{m\in[0,1]}\left(m^2\frac{\pi}{r}+\frac{1}{C^*\varepsilon}(1-m)^2 \right)=\frac{\pi}{r}\frac{1}{1+\pi C^*\varepsilon/r},
\end{equation}
with $c_0,C^*$ constants defined in Lemma 3.1,3 and Lemma 3.1.4 in \cite{CollianderJerrard1999GLvorticesSE}.
Since $r>r_\varepsilon>r_0=\delta_\varepsilon^{1/4}/2$, $\lambda^\varepsilon(r)<c_0/\varepsilon$ for $\varepsilon$ small enough. Then \eqref{eq:Eue B1-Bre} gives
\begin{align}
	(1+g)\int_{B_1\setminus B_{r_\varepsilon}(\bvec{0})}e_0^\varepsilon(u^\varepsilon,v^\varepsilon)\ud\vx\ge&\int_{r\in T_\varepsilon}\frac{\pi}{r}\frac{1}{1+\pi C^*\varepsilon/r}\ud r\ge\int_{r_\varepsilon+C\delta_\varepsilon}^1\frac{\pi}{r}\frac{1}{1+\pi C^*\varepsilon/r}\ud r\nn =&\pi\log\frac{1+\pi C^*\varepsilon}{r_\varepsilon+C\delta_\varepsilon+\pi C^*\varepsilon}\ge \pi\log\frac{1}{r_\varepsilon}- C\delta_\varepsilon^{3/4},
\end{align}
which together with \eqref{eq:Etue B1} implies
\begin{equation}
	\int_{B_1(\bvec{0})}e_0^\varepsilon(u^\varepsilon,v^\varepsilon)\ud\vx\ge \int_{B_1(\bvec{0})}e^\varepsilon_0(\tilde{u}^\varepsilon,\tilde{v}^\varepsilon)\ud \vx-C\tilde{\delta}_\varepsilon.
\end{equation}
Noting \eqref{eq:def of gammag} and Lemma {lem:existence of gammag}, we have
\begin{equation}
	\int_{B_1}e^\varepsilon_0(u^\varepsilon,v^\varepsilon)\ud\vx\ge \frac{\pi}{1+g}\log\frac{1}{\varepsilon\sqrt{1+g}}+\gamma_g-O(\varepsilon^2)-C\tilde{\delta}_\varepsilon,
\end{equation}
which complete the proof.
\end{proof}
\section{On the ODE related to the stationary state}\label{sec:ODE}
In this section, we will prove the existence of $\gamma_g$.

\begin{Lem}\label{lem:existence of gammag}
There exists $\gamma_g\in\R$ such that 
\begin{equation}\label{eq:gammag=Ie-loge+e2}
	\gamma_g:=\inf_{(u,v)\in A_g(B_1(\bvec{0}))}\int_{B_1(\bvec{0})}e^\varepsilon_0(u,v)\ud\vx-\frac{\pi}{1+g}\log\frac{1}{\varepsilon\sqrt{1+g}}+O(\varepsilon^2),
\end{equation}
with
\[
	A_g(B_1(\bvec{0}))=\left\{(u,v)\in H^1(B_1(\bvec{0}))\times H^1(B_1(\bvec{0}))\left|(u(\vx),v(\vx))=\frac{1}{\sqrt{1+g}}\left(\frac{z}{|z|},1\right)\  \text{for}\ \vx\in \p B_1(\bvec{0})\right.\right\}.
\]
\end{Lem}
In fact, $\inf_{(u,v)\in A_g(B_1(\bvec{0}))}\int_{B_1(\bvec{0})}e^\varepsilon_0(u,v)\ud\vx=\int_{B_1(\bvec{0})}e^\varepsilon_0(u_m,v_m)\ud\vx$ where $u_m(\vx)=f_1(|\vx|)e^{i\theta(\vx)},v_m(\vx)=f_2(|\vx|)$ with $e^{i\theta(\vx)=(x+iy)/|\vx|}$ and $f_1,f_2$ solving
\begin{equation}\label{eq:ODE of norm e}
	\begin{cases}
	-f_1''(r)-\frac{1}{r}f_1'(r)+\frac{1}{r^2}f_1(r)+\frac{1}{\varepsilon^2}(f_1(r)^2+g f_2(r)^2-1)f_1(r)=0,&0<r<1,\\
	-f_2''(r)-\frac{1}{r}f_2'(r)+\frac{1}{\varepsilon^2}(f_2(r)^2+g f_1(r)^2-1)f_2(r)=0,&0<r<1,\\
	f_1(0)=0=f_2(0),\quad f_1(1)=\frac{1}{\sqrt{1+g}}=f_2(1).
	\end{cases}
\end{equation}
Moreover, 
\begin{align}
	&\inf_{(u,v)\in A_g(B_1(\bvec{0}))}\int_{B_1(\bvec{0})}e^\varepsilon_0(u,v)\ud\vx\nn &\quad=\pi\int_0^1\left( r|f_1'|^2+r|f_2'|^2+\frac{f_1^2}{r}\right.\nn &\qquad\qquad+\left.\frac{r}{2\varepsilon^2}\left( (f_1^2-\frac{1}{1+g})^2+2g(f_1^2-\frac{1}{1+g})(f_2^2-\frac{1}{1+g})+(f_2^2-\frac{1}{1+g})^2\right)\right)\ud r.
\end{align}
Hence, we only need to verify the existence and convergence of the solution of \eqref{eq:ODE of norm e}. 
On the other hand, via a rescaling of the solution $(f_1,f_2)$ of \eqref{eq:ODE of norm e}:
\begin{equation}
	\tilde{f}_1(r)=f_1(r\varepsilon),\tilde{f}_2(r)=f_2(r\varepsilon),r\in [0,1/\varepsilon],
\end{equation}
we will obtain a solution $(\tilde{f}_1,\tilde{f}_2)$ of the solution 
\begin{equation}\label{eq:ODE of norm R}
	\begin{cases}
	-f_1''(r)-\frac{1}{r}f_1'(r)+\frac{1}{r^2}f_1(r)+(f_1(r)^2+g f_2(r)^2-1)f_1(r)=0,&0<r<R,\\
	-f_2''(r)-\frac{1}{r}f_2'(r)+(f_2(r)^2+g f_1(r)^2-1)f_2(r)=0,&0<r<R,\\
	f_1(0)=0=f_2'(0),\quad f_1(R)=\frac{1}{\sqrt{1+g}}=f_2(R),
	\end{cases}
\end{equation}
with $R=1/\varepsilon$. Hence, we only need to study the existence and convergence of the solution of \eqref{eq:ODE of norm R}, which are the main consideration of this section.

\begin{Lem}[Existence of the solution of ODE]\label{lem:existence of ODER}
For any $0<R<\infty$, there exists $(f_1,f_2)$ solving \eqref{eq:ODE of norm R} without the condition $f_1(0)=0=f'_2(0)$. Moreover, $0\le f_1\le 1/\sqrt{1+g}\le f_2\le 1$.
\end{Lem}
\begin{proof}
For any $(h_1,h_2)\in M_R:=\{(h_1,h_2):(0,R)\to \R^2|\sqrt{r}h_1',h_1/\sqrt{r},h_2'\in L^2(0,R),h_1(R),h_2(R)=1/\sqrt{1+g} \}$,
\begin{equation}
	E_R(h_1,h_2):=\int_0^R(e_k(h_1,h_2)+e_p(h_1,h_2))\ud r,
\end{equation}
where 
\begin{equation}
	e_k(h_1,h_2)=r|h_1'|^2+r|h_2'|^2+\frac{h_1^2}{r},e_p(h_1,h_2)=\frac{r}{2}(h_1^2-\frac{1}{1+g})^2+gr(h_1^2-\frac{1}{1+g})(h_2^2-\frac{1}{1+g})+\frac{r}{2}(h_2^2-\frac{1}{1+g})^2.
\end{equation}
To prove the existence of the solution of \eqref{eq:ODE of norm R}, we only need to prove the existence of minimizer of $E_R(h_1,h_2)$ in $M_R$.

% We first prove that if the minimizer $(f_1,f_2)$ exists, then 
% \begin{equation}
% 	(f_1,f_2)\in \tilde{M}_R:=\{(f_1,f_2)\in M_R|0\le f_1\le 1/\sqrt{1+g}\le f_2\le 1\}.
% \end{equation}
For any $(h_1,h_2)\in M_R$, we define 
\begin{equation}\label{eq:def of th}
	\tilde{h}_1=\min\left\{|h_1|,\frac{1}{\sqrt{1+g}}\right\},\tilde{h}_2=\begin{cases}1/\sqrt{1+g},&|h_2|<1/\sqrt{1+g},\\ |h_2|,&1/\sqrt{1+g}\le |h_2|\le 1,\\1,&|h_2|>1. \end{cases}
\end{equation}
\eqref{eq:def of th} immediately gives
\begin{equation}\label{eq:ekth<ekh}
	e_k(\tilde{h}_1,\tilde{h}_2)\le e_k(h_1,h_2).
\end{equation}
And we can consider $e_p(\tilde{h}_1,\tilde{h}_2)$ case by case to obtain 
\begin{equation}\label{eq:epth < eph}
	e_p(\tilde{h}_1,\tilde{h}_2)\le e_p(h_1,h_2).
\end{equation}
For example, for the case $|h_1|<1/\sqrt{1+g},|h_2|>1$, 
\begin{equation}
	e_p(h_1,h_2)-e_p(\tilde{h}_1,\tilde{h}_2)=r(h_2^2-1)\left(h_2^2+1+2gh_1^2-\frac{4g}{1+g}\right)\ge r(h_2^2-1)\left(2-\frac{4g}{1+g}\right)=r(h_2^2-1)\frac{2-2g}{1+g}\ge 0.
\end{equation}
Combining \eqref{eq:ekth<ekh} and \eqref{eq:epth < eph}, we have
\begin{equation}
	E_R(\tilde{h}_1,\tilde{h}_2)\le E_R(h_1,h_2),
\end{equation}
which implies
\begin{equation}
	\inf_{(h_1,h_2)\in M_R}E_R(h_1,h_2)=\inf_{(h_1,h_2)\in \tilde{M}_R}E_R(h_1,h_2),
\end{equation}
with
\[
	\tilde{M}_R=\{(f_1,f_2)\in M_R|0\le f_1\le 1/\sqrt{1+g}\le f_2\le 1\}.
\]
Since $E_R(h_1,h_2)\ge 0$, we can find a sequence $(h_1^k,h_2^k)\in \tilde{M}_R$ such that 
\begin{equation}
	\lim_{k\to\infty}E_R(h_1^k,h_2^k)=\inf_{(h_1,h_2)\in \tilde{M}_R}E_R(h_1,h_2)\ge 0.
\end{equation}
Then \[
	\|\sqrt{r}(h_1^k)'\|_{L^2(0,R)}^2+\|\sqrt{r}(h_2^k)'\|_{L^2(0,R)}^2+\|h_1^k/\sqrt{r}\|_{L^2(0,R)}^2\le \sup_{k\ge 1} E_R(h_1^k,h_2^k) \le C,
\]
which implies that there exist $H_1,H_2,H_3\in L^2(0,R)$ such that (up to a subsequence)
\begin{equation}
	\sqrt{r}(h_1^k)'\wto H_1,\sqrt{r}(h_2^k)'\wto H_2,h_1^k/\sqrt{r}\wto H_3,\quad \text{in}\ L^2(0,R).
\end{equation}
We define 
\begin{equation}
	f_1=\sqrt{r}H_3,f_2=\frac{1}{\sqrt{1+g}}-\int_r^R\frac{H_2(s)}{\sqrt{s}}\ud s.
\end{equation}
Then we have $H_1=\sqrt{r}f_1'$ and hence $h_1^k\wto f_1$ in $H^1([R_0,R])\hookrightarrow\hookrightarrow C^{0.5}(R_0,R)$ for any $0<R_0<R$. As a result $h_1^k\to f_1$ in $C^{0.5}_{loc}(0,R)$ and $f_1(R)=1/\sqrt{1+g},1\le f_1\le 1/\sqrt{1+g}$. On the other hand, 
\begin{equation}
	h_2^k(r)=\sqrt{r}H_3,f_2=\frac{1}{\sqrt{1+g}}-\int_r^R\frac{\sqrt{s}(h_2^k)'(s)}{\sqrt{s}}\ud s.
\end{equation}
Letting $k\to \infty$ and noting the definition of $f_2$, we have $\lim_{k\to\infty}h_2^k(r)=f_2(r)$ and hence $1/\sqrt{1+g}\le f_2\le 1,f_2(R)=1/\sqrt{1+g}$. As a result $(f_1,f_2)\in M_R$. Noting 
\begin{equation}
	\sqrt{r}(h_1^k)'\wto \sqrt{r}f_1',\sqrt{r}(h_2^k)'\wto \sqrt{r}f_2',h_1^k/\sqrt{r}\wto f_1/\sqrt{r},\quad \text{in}\ L^2(0,R),
\end{equation}
and 
\begin{equation}
	h_1^k(r)\to f_1(r),h_2^k(r)\to f_2(r), \ \text{as}\ k\to \infty \ \text{for a.e.}\ r\in [0,R],
\end{equation}
we finally obtain 
\begin{equation}
	\inf_{(h_1,h_2)\in \tilde{M}_R}E_R(h_1,h_2)\le E_R(f_1,f_2)\le \lim_{k\to \infty}E_R(h_1^k,h_2^k)=\inf_{(h_1,h_2)\in \tilde{M}_R}E_R(h_1,h_2),
\end{equation}
which complete the proof.
\end{proof}

\begin{Lem}\label{lem:f10=0,f2'0=0 R}
The solution obtained in Lemma \ref{lem:existence of ODER} satisfies
$f_1(r)\le Cr$ for $r$ small enough and $f_2'(0)=0$, i.e. it's a solution to \eqref{eq:ODE of norm R}.
\end{Lem}

\begin{proof}
Since $0\le f_1\le 1/\sqrt{1+g}\le f_2\le 1$, we have
\begin{equation}
	-f_1''(r)-\frac{1}{r}f_1'(r)+\frac{1}{r^2}f_1(r)=-(f_1^2+gf_2^2-1)f_1\le Cf_1,
\end{equation}
i.e. $Lf_1\le 0$, with
\begin{equation}
	Lh:=-h''-\frac{h'}{r}+\left(\frac{1}{r^2}-C \right)h.
\end{equation}
Noting that there are two functions $w_1(r)=r+o(r)$ and $w_2(r)=r^{-1}+0(r^-1)$ satisfying $Lw_1=Lw_2=0$,we can obtain there are some $\mu>0$ such that $f_1(r)\le \mu r$ via repeating Step 1 in the proof of Proposition 2.2 in \cite{IgnatNguyenSlastikovValeriy2014UniquenessODEGinzburgLandau}. 

For the proof of $f_2'(0)=0$, we consider $h(r)=(f_2^2+gf_1^2-1)f_2$, then we have 
\begin{equation}
	(rf_2)'=rf_2'+f_2=rh(r)\le Mr \ \text{in}\ (0,R).
\end{equation}
Then integrating on $(R_0,r)$ we have
\begin{equation}
	rf_2'(r)=R_0f_2'(R_0)+\int_{R_0}^2sf(s)\ud s,
\end{equation}
which implies $\lim_{\varepsilon\to 0^+}rf_2(r)$ exists and we denote it by $c_0$. If $c_0\ne 0$, without loss of generality, we assume $c_0>0$. Then there exists $\delta_0>0$ s.t. $rf_2'(r)\ge c_0/2$ for any $r\in (0,\delta_0)$. Then we have
\begin{equation}
	f_2(r)=f_2(\delta_0)+\int_{\delta_0}^rf_2'(s)\ud s\le f_2(\delta_0)-\frac{c_0}{2}\log\frac{\delta_0}{r},
\end{equation}
which contradicts to $f_2(r)\ge 1/\sqrt{1+g}$. As a result, $c_0=0$ and
\begin{equation}
	rf_2'(r)=\int_0^rsf(s)\ud s.
\end{equation}
Since $h(s)\le M$, we have
\begin{equation}
	|f_2'(r)|\le \frac{1}{r}\int_0^rsM\ud s\le \frac{Mr}{2},
\end{equation}
which completes the proof.
\end{proof}

\begin{Lem}\label{lem:fjR-s=r-2}
The solution obtained in Lemma \ref{lem:existence of ODER} satisfies 
\begin{equation}
	f_1(r)\ge \frac{1}{\sqrt{1+g}}-\frac{C}{r^2}, \quad f_2(r)\le \frac{1}{\sqrt{1+g}}+\frac{C}{r^2}
\end{equation}
for $R,r$ large enough.
\end{Lem}
\begin{proof}
% We consider 
% \begin{equation}
% 	\varphi_1(\tau)=\frac{1}{\sqrt{1+g}}- f_1(1/\tau),\varphi_2(\tau)=f_2(1/\tau)-\frac{1}{\sqrt{1+g}},\quad \tau\ge 1/R.
% \end{equation}
% Then $\varphi_1,\varphi_2$ satisfy
% \begin{equation}
% \begin{cases}
% 	-\varphi_1''(\tau)-\frac{1}{\tau}\varphi_1'(\tau)+\frac{1}{\tau^2}\varphi_1(\tau)=\frac{1}{\tau^2\sqrt{1+g}}+\frac{1}{\tau^4}(\varphi_1^2+g\varphi_2^2-\frac{2\varphi_1}{\sqrt{1+g}}+\frac{2g\varphi_2}{\sqrt{1+g}})(\frac{1}{\sqrt{1+g}}-\varphi_1),\\
% 	-\varphi_2''(\tau)-\frac{1}{\tau}\varphi_2'(\tau)=-\frac{1}{\tau^4}(\varphi_2^2+g\varphi_1^2-\frac{2g\varphi_1}{\sqrt{1+g}}+\frac{2\varphi_2}{\sqrt{1+g}})(\frac{1}{\sqrt{1+g}}+\varphi_2).
% 	\end{cases}
% \end{equation}
% Step 1: We will prove that $f_1^2+f_2^2\le 2/(1+g)$.

Let $S=f_1^2+f_2^2$. Then $S$ satisfies
\begin{equation}
	-S''-\frac{1}{r}S'+\frac{S}{r^2}+((1+g)S-2)S=-(1-g)(f_1^2-f_2^2)^2-2((f_1')^2+(f_2')^2)\le 0.
\end{equation}
Noting $S(0)=f_1^2(0)+f_2^2(0)\le 1\le 2/(1+g)$ and $S(R)=2/(1+g)$, we can prove $S(r)\le 2/(1+g)$ via the maximum principle. Equivalently, $f_2^2-1/(1+g)\le 1/(1+g)-f_1^2$, which implies
\begin{equation}
	0=-f_1''-\frac{1}{r}f_1'+\frac{1}{r^2}f_1+(f_1^2+gf_2^2-1)f_1\le -f_1''-\frac{1}{r}f_1'+\frac{1}{r^2}f_1+(1-g)(f_1^2-\frac{1}{1+g})f_1,
\end{equation}
i.e. $f_1$ is a supersolution of $\tilde{L}f=0$ where
\begin{equation}
	\tilde{L}f:=-f''-\frac{1}{r}f'+\frac{1}{r^2}f+(1-g)(f^2-\frac{1}{1+g})f.
\end{equation}
Via Theorem 1.1 in \cite{IgnatNguyenSlastikovValeriy2014UniquenessODEGinzburgLandau}, there are $\tilde{f}^R$ and $\tilde{f}$ satisfying 
\begin{equation}
	\tilde{L}\tilde{f}^R=0, r\in (0,R) ,\quad \tilde{f}^R(0)=0,\tilde{f}^R(R)=\frac{1}{\sqrt{1+g}},
\end{equation}
and 
\begin{equation}
	\tilde{L}\tilde{f}=0, r\in (0,+\infty) ,\quad \tilde{f}(0)=0,\tilde{f}(+\infty)=\frac{1}{\sqrt{1+g}}.
\end{equation}
Moreover, $\tilde{f}(R)\le 1/\sqrt{1+g}=\tilde{f}^R(R)$. Hence, noting Proposition 3.5 in \cite{IgnatNguyenSlastikovValeriy2014UniquenessODEGinzburgLandau} and Proposition 2.5, there is some constant $C$ and $R_0$ such that if $R>R_0$ and $r>R_0$, 
\begin{equation}
	f_1(r)\ge\tilde{f}^R(r)\ge\tilde{f}(r)\ge \frac{1}{\sqrt{1+g}}-\frac{C}{r^2}.
\end{equation}
Then, noting $f_2^2-1/(1+g)\le 1/(1+g)-f_1^2$ and $f_1\le \sqrt{1}{1+g}\le f_2$, we have
\begin{equation}
	f_2-\frac{1}{\sqrt{1+g}}\le \frac{1}{\sqrt{1+g}}-f_1\le \frac{C}{r^2},
\end{equation}
which completes the proof.
% Step 2: Noting 
\end{proof}
% \begin{Lem}
% The solution obtained in Lemma \ref{lem:existence of ODER} satisfies 
% \begin{equation}
% 	f_1(r)=\frac{1}{\sqrt{1+g}}-\frac{\alpha}{r^2}+O(r^{-4}),\quad f_2(r)=\frac{1}{\sqrt{1+g}}+\frac{\beta}{r^2}+O(r^-4),\quad r\to \infty.
% \end{equation}
% with $\alpha=\frac{\sqrt{1+g}}{2(1-g^2)},\beta=\frac{g\sqrt{1+g}}{2(1-g^2)}$
% for $R,r$ large enough.
% \end{Lem}
% \begin{proof}
% We consider 
% \begin{equation}\label{eq:def of f11/tau}
% 	\varphi_1(\tau)=\frac{1}{\sqrt{1+g}}-\alpha \tau^2- f_1(1/\tau),\varphi_2(\tau)=f_2(1/\tau)-\frac{1}{\sqrt{1+g}}-\beta \tau^2,\quad \tau\ge 1/R.
% \end{equation}
% Substituting \eqref{eq:def of f11/tau} to \eqref{eq:ODE of norm R} we obtain
% Then $\varphi_1,\varphi_2$ satisfy
% \begin{equation}
% \begin{cases}
% 	-\varphi_1''(\tau)-\frac{1}{\tau}\varphi_1'(\tau)+\frac{1}{\tau^2}\varphi_1(\tau)+3\alpha+\frac{1}{\tau^4}\frac{2}{1+g}() =\frac{1}{\tau^2\sqrt{1+g}}+\frac{1}{\tau^4}(\varphi_1^2+g\varphi_2^2-\frac{2\varphi_1}{\sqrt{1+g}}+\frac{2g\varphi_2}{\sqrt{1+g}})(\frac{1}{\sqrt{1+g}}-\varphi_1),\\
% 	-\varphi_2''(\tau)-\frac{1}{\tau}\varphi_2'(\tau)=\frac{1}{\tau^4}(\varphi_2^2+g\varphi_1^2-\frac{2g\varphi_1}{\sqrt{1+g}}+\frac{2\varphi_2}{\sqrt{1+g}})(\frac{1}{\sqrt{1+g}}+\varphi_2).
% 	\end{cases}
% \end{equation}
% \end{proof}

\begin{Lem}\label{lem:existence of solution of norm inf}
Denote $(f_1^R,f_2^R)$ to be the solution of \eqref{eq:ODE of norm R}, then as $R\to \infty$, (up to a subsequence)
\begin{equation}
	(f_1^R,f_2^R)\to (f_1,f_2) \in C_{loc}^2\left((0,\infty)\right)
\end{equation}
and $(f_1,f_2)$ solves 
\begin{equation}\label{eq:ODE of norm infty}
	\begin{cases}
	-f_1''(r)-\frac{1}{r}f_1'(r)+\frac{1}{r^2}f_1(r)+(f_1(r)^2+g f_2(r)^2-1)f_1(r)=0,&0<r<\infty,\\
	-f_2''(r)-\frac{1}{r}f_2'(r)+(f_2(r)^2+g f_1(r)^2-1)f_2(r)=0,&0<r<\infty,\\
	f_1(0)=0=f_2'(0),\quad f_1(\infty)=\frac{1}{\sqrt{1+g}}=f_2(\infty).
	\end{cases}
\end{equation}
Moreover, $f_1(r)=1/\sqrt{1+g}+O(r^-2),f_2(r)=1/\sqrt{1+g}+O(r^-2)$.
\end{Lem}
\begin{proof}
Any compact interval $(r_1,r_2)\subset(0,+\infty)$ would be contained in $(0,R)$ for $R$ large enough.

Then $(f_1^R,f_2^R)$ is a minimizer of $\int_{r_1}^{r_2}e_R(h_1,h_2)\ud r$ in $M_{R,r_1,r_2}:=\{(h_1,h_2)\in M_R:h_j(r_k)=f_j^R(r_k)\ \text{for}\ 1\le j,k\le 2\}$. we take $h_j(r)=(f_j^R(r_2)-f_j^R(r_1))(r-r_1)/(r_2-r_1)+f_j^R(r_1)$. Then  
\begin{equation}
	\int_{r_1}^{r_2}e_R(f_1^R,f_2^R)\ud r\le \int_{r_1}^{r_2}e_R(h_1,h_2)\ud r\le \frac{1}{1+g}\left(\log\frac{r_2}{r_1}+\frac{1}{2}\left(\frac{r_2-r_1}{r_2+r_1}+r_2^2-r_1^2 \right) \right),
\end{equation}
i.e. there is a uniform upper bound of $\int_{r_1}^{r_2}e_R(h_1,h_2)\ud r$ in $M_{R,r_1,r_2}$ for any $R>R_0$. 
Then we will have $\| f_1^R \|_{H^4(I)}+\| f_2^R\|_{H^4(I)}$ is uniformly bounded via the equation \eqref{eq:ODE of norm R}, and hence $\| f_1^R \|_{C^3(I)}+\| f_2^R\|_{C^3(I)}$ is uniformly bounded. Then  $f_1^R\to f_1,f_2\to f_2$ in $C^2(I)$ (up to a subsequence) for some $f_j\in C^2(I)$. Noting that $I$ is arbitrary, we have $f_1^R\to f_1,f_2\to f_2$ in $C^2_{loc}\left((0,\infty)\right)$ (up to a subsequence) for some $f_j\in C^2\left((0,\infty)\right)$.
Then $(f_1,f_2)$ will satisfy \eqref{eq:ODE of norm infty} in $(0,\infty)$. Similar to Lemma \ref{lem:f10=0,f2'0=0 R}, $f_1(0)=0,f_2'(0)=0$. And Lemma \ref{lem:fjR-s=r-2} implies $|f_j^R(r)-1/\sqrt{1+g}|\le Cr^{-2}$ for $r>R_0,R>R_0$. Hence $|f_j(r)-1/\sqrt{1+g}|\le Cr^{-2}$ for $r>R_0$, which implies $f_j(r)=1/\sqrt{1+g}+O(r^{-2})$.
\end{proof}

\begin{Lem}\label{lem:f=s+Cr2+Cr4}
For the solution of \eqref{eq:ODE of norm infty}, we have
\begin{equation}\label{eq:f=s+Cr2+Cr4}
	f_1(r)=\frac{1}{\sqrt{1+g}}-\frac{\alpha}{r^2}+O(r^{-4}),\quad f_2(r)=\frac{1}{\sqrt{1+g}}+\frac{\beta}{r^2}+O(r^{-4}),\quad r\to \infty,
\end{equation}
with $\alpha=\frac{\sqrt{1+g}}{2(1-g^2)},\beta=\frac{g\sqrt{1+g}}{2(1-g^2)}$.
\end{Lem}
\begin{proof}
We denote $s_+=1/\sqrt{1+g}$ and consider 
\begin{equation}\label{eq:def of f11/tau}
	\varphi_1(\tau)=f_1(1/\tau)-\frac{1}{\sqrt{1+g}}+\alpha \tau^2,\varphi_2(\tau)=f_2(1/\tau)-\frac{1}{\sqrt{1+g}}-\beta \tau^2,\quad \tau>0.
\end{equation}
Substituting \eqref{eq:def of f11/tau} to \eqref{eq:ODE of norm R} we obtain
% \begin{equation}
% \begin{cases}
% 	-\varphi_1''(\tau)-\frac{1}{\tau}\varphi_1'(\tau)+\frac{1}{\tau^2}\varphi_1(\tau)+3\alpha+\frac{1}{\tau^4}\frac{2}{1+g}() =\frac{1}{\tau^2\sqrt{1+g}}+\frac{1}{\tau^4}(\varphi_1^2+g\varphi_2^2-\frac{2\varphi_1}{\sqrt{1+g}}+\frac{2g\varphi_2}{\sqrt{1+g}})(\frac{1}{\sqrt{1+g}}-\varphi_1),\\
% 	-\varphi_2''(\tau)-\frac{1}{\tau}\varphi_2'(\tau)=\frac{1}{\tau^4}(\varphi_2^2+g\varphi_1^2-\frac{2g\varphi_1}{\sqrt{1+g}}+\frac{2\varphi_2}{\sqrt{1+g}})(\frac{1}{\sqrt{1+g}}+\varphi_2).
% 	\end{cases}
% \end{equation}
% \begin{equation}
% 	\begin{cases}
% 	-\varphi_1''-\frac{1}{\tau}\varphi_1'+\frac{1}{\tau^2}\varphi_1+\left(\frac{1}{\tau^4}(2s_+^2+3s_+\varphi_1+2s_+g\varphi_2+\varphi_1^2+g\varphi_2^2)+\frac{1}{\tau^2}\left(-4s_+\alpha+2g\beta\varphi_2-3\alpha\varphi_1\right)+(3\alpha^2+g\beta^2)\right)\varphi_1\nn 
% 	\quad +\left(\frac{1}{\tau^4}(2s_+^2g+s_+g\varphi_2)+\frac{1}{\tau^2}(2s_+g(\beta-\alpha)-2g\alpha \beta)\varphi_2\right)=-3\alpha-s_+(3\alpha^2-2g\alpha\beta+g\beta^2)+\alpha(\alpha^2+g\beta^2)\tau^2 ,&\tau>0,
% 	\end{cases}
% \end{equation}
\begin{equation}
	\begin{cases}
	-\varphi_1''-\frac{1}{\tau}\varphi_1'+\frac{2s_+^2}{\tau^4}\varphi_1+\frac{2s_+^2g}{\tau^4}\varphi_2=R_1(\tau),&\tau>0,\\
	-\varphi_2''-\frac{1}{\tau}\varphi_2'+\frac{2s_+^2}{\tau^4}\varphi_2+\frac{2s_+^2g}{\tau^4}\varphi_1=R_2(\tau),&\tau>0,
	\end{cases}
\end{equation}
where 
\begin{equation}
\begin{cases}
 	R_1(\tau):=-F_3(-\alpha,\beta,\tau)-3\alpha-\frac{1}{\tau^2}\varphi_1-F_1(-\alpha,\beta,\varphi_1,\varphi_2,\tau)\frac{\varphi_1}{\tau^4}-F_2(-\alpha,\beta,\varphi_1,\varphi_2,\tau)\frac{\varphi_2}{\tau^4} ,\\
 	R_2(\tau):=-F_3(\beta,-\alpha,\tau)+4\beta-F_1(\beta,-\alpha,\varphi_2,\varphi_1,\tau)\frac{\varphi_2}{\tau^4}-F_2(\beta,-\alpha,\varphi_2,\varphi_1,\tau)\frac{\varphi_1}{\tau^4},
 	 \end{cases}
 \end{equation} and for $x_1,x_2,y_1,y_2,\tau\in \R$,
\begin{equation}
\begin{cases}
	F_1(x_1,x_2,y_1,y_2,\tau)=3s_+y_1+2s_+gy_2+y_1^2+gy_2^2+\tau^2(4s_+x_1+2gx_2y_2+3x_1y_1)+\tau^4(3x_1^2+gx_2^2),\\
	F_2(x_1,x_2,y_1,y_2,\tau)=s_+gy_2+\tau^2(2s_+g(x_1+x_2)+gx_1y_2)+2gx_1x_2\tau^4,\\
	F_3(x_1,x_2,\tau)=s_+(3x_1^2+2gx_1x_2+gx_2^2)+x_1(x_1^2+gx_2^2)\tau^2.
	\end{cases}
\end{equation}
Noting that $|F_j(x_1,x_2,y_1,y_2,\tau)|\le C(|y_1|+|y_2|+\tau^2),j=1,2,|F_3(x_1,x_2,\tau)|\le C(x_1^2+x_2^2+\tau^2)$ for $|y_1|,|y_2|,|x_1|,|x_2|,\tau^2\le |\alpha|$ and that $|\varphi_1(\tau)|+|\varphi_2(\tau)|\le C\tau^2$ when $\tau$ is small enough. Hence we can find $\delta>0$ such that for any $0<\tau<\delta$, 
\begin{equation}
	|R_1(\tau)|\le C_*, |R_2(\tau)|\le C_*,16\tau^2<1.
\end{equation}
We define
\begin{equation}\begin{cases}
	L_1(h_1,h_2)=-h_1''-\frac{1}{\tau}h_1'+\frac{2s_+^2}{\tau^4}h_1+\frac{2s_+^2g}{\tau^4}h_2,\\
	L_2(h_1,h_2)=-h_2''-\frac{1}{\tau}h_2'+\frac{2s_+^2}{\tau^4}h_2+\frac{2s_+^2g}{\tau^4}h_1,\\
	D=\max\left\{C,\frac{|\varphi_1(\delta)|}{\delta^4},\frac{|\varphi_2(\delta)|}{\delta^4} \right\}.
	\end{cases}
\end{equation}
Then, we have
\begin{equation}
	L_k(\phi_1,\phi_2)\ge C,\quad \phi_1(\delta)\ge|\varphi_1(\delta)|,\quad \phi_2(\delta)\ge|\varphi_2(\delta)|,
\end{equation}
where $\phi_1(\tau)=D\tau^4,\phi_2(\tau)=D\tau^4$. If we define $V_j=\phi_j-\varphi_j$, 
\begin{equation}
	L_k(V_1,V_2)\ge 0, V_j(0)\ge 0,V_j(\delta)\ge 0.
\end{equation}
Summing up $L_1(V_1,V_2),L_2(V_1,V_2)$ we obtain that
\begin{equation}
	-(V_1+V_2)''-\frac{1}{\tau}(V_1+V_2)'+\frac{2s_+^2(1+g)}{\tau^4}(V_1+V_2)\ge0, V_1(0)+V_2(0)\ge 0, V_1(\delta)+V_2(\delta)\ge0,
\end{equation}
which implies $V_1(\tau)+V_2(\tau)\ge0$ in $(0,\delta)$, i.e. $\varphi_1(\tau)+\varphi_2(\tau)\le 2D\tau^4$. Similarly, we take $V_1=\phi_1-\varphi_1,V_2=\phi_2+\varphi_2$ to obtain $\varphi_1(\tau)-\varphi_2(\tau)\le 2D\tau^4$, taking $V_1=\phi_1+\varphi_1,V_2=\phi_2+\varphi_2$ gives $-\varphi_1(\tau)-\varphi_2(\tau)\le 2D\tau^4$, and we have $-\varphi_1(\tau)+\varphi_2(\tau)\le 2D\tau^4$ by taking $V_1=\phi_1+\varphi_1,V_2=\phi_2-\varphi_2$. Then we have $|\varphi_1(\tau)|+|\varphi_2(\tau)|\le 2D\tau^4$ in $(0,\delta)$, which complete the proof by \eqref{eq:f=s+Cr2+Cr4}.
\end{proof}

\begin{Lem}\label{lem:f'=Cr3}
For the solution of \eqref{eq:ODE of norm infty}, we have
\begin{equation}\label{eq:f'=Cr3}
	f_1'(r)=O(r^{-3}),\quad f_2'(r)=O(r^{-3}),\quad r\to \infty.
\end{equation}
\end{Lem}
\begin{proof}
Let 
We denote $s_+=1/\sqrt{1+g}$ and consider 
\begin{equation}\label{eq:def of f1r}
	\varphi_1(r)=f_1(r)-\frac{1}{\sqrt{1+g}}+\frac{\alpha}{r^2}, \varphi_2(r)=f_2(r)-\frac{1}{\sqrt{1+g}}-\frac{\beta}{r^2},\quad r>0.
\end{equation}
Lemma \ref{lem:f=s+Cr2+Cr4} implies that there exist $C_1>0,R_0>0$ such that $|\varphi_1(r)|=C_1r^{-4},|\varphi_2(r)|=C_1r^{-4}$ for  $r>R_0>0$. Then we have
\begin{equation}
	|-\varphi_1''-\frac{1}{r}\varphi_1'=Cr^{-4}, -\varphi_2''-\frac{1}{r}\varphi_2'|\le Cr^{-4},\quad r>R_0
\end{equation}
i.e.
\begin{equation}
	|(r\varphi_j')'|\le Cr^{-3},\quad j=1,2, r>R_0
\end{equation}
For any $r_1>r_2>R_0$, we have
\begin{equation}\label{eq:rf'r-Rf'R}
	-\frac{C}{2r_2^2}+\frac{C}{2r_1^2}\le r_1\varphi'_j(r_1)-r_2\varphi_j'(r_2)\le\frac{C}{2r_2^2}-\frac{C}{2r_1^2}.
\end{equation}
Then there must be $r\varphi_j'(r)\to 0$ as $r\to \infty$. Otherwise, there is some $c>0$ and a sequence $r_k\to\infty$ such that
\begin{equation}
	|r_k\varphi_j'(r_k)|\ge c.
\end{equation}
Without loss of generality, we can assume
\begin{equation}
	r_k\varphi_j'(r_k)\ge c.
\end{equation}
Then for $r>R_0$, we can always find $r_k>r$ and hence we have
\begin{equation}
	r\varphi_j'(r)\ge r_k\varphi_j'(r_k)-\frac{C}{2r^2}+\frac{C}{2r_k^2}\ge c-\frac{C}{2r_1^2}+\frac{C}{2r_k^2}.
\end{equation}
Then for any $r\ge R_1:=\max\{R_0,\sqrt{C/c}\}$, we have $r\varphi_j'(r)\ge c/2$, which immediately gives
\begin{equation}
	\varphi_j(r)\ge \varphi_j(R_1)+\int_{R_1}^r\varphi_j'(r)\ud r\ge \varphi_j(R_1)+\frac{c}{2}\log\frac{r}{R_1},
\end{equation}
contradicting to $\varphi_j(r)=O(r^{-4})$. Letting $r_1\to \infty$ in \eqref{eq:rf'r-Rf'R}, we obtain
\begin{equation}
	-\frac{C}{2r_2^2}\le -r_2\varphi_j'(r_2)\le\frac{C}{2r_2^2}, \quad r_2 >R_0.
\end{equation}
Noting \eqref{eq:def of f1r}, we have
\begin{equation}
	f_1'(r)=\frac{2\alpha}{r^3}+\varphi_1'(r)=O(r^{-3}),f_1'(r)=-\frac{2\beta}{r^3}+\varphi_2'(r)=O(r^{-3}).
\end{equation}
\end{proof}

With above lemmas, we can proof Lemma \ref{lem:existence of gammag}.
\begin{proof}
We denote
\begin{equation}
	I(\varepsilon,1)=\inf_{(u,v)\in A_g(B_1(\bvec{0}))}\int_{B_1(\bvec{0})}e^\varepsilon_0(u,v)\ud\vx
\end{equation}
with 
\[
	A_g(B_1(\bvec{0}))=\left\{(u,v)\in H^1(B_1(\bvec{0}))\times H^1(B_1(\bvec{0}))\left|(u(\vx),v(\vx))=\frac{1}{\sqrt{1+g}}\left(\frac{z}{|z|},1\right)\  \text{for}\ \vx\in \p B_1(\bvec{0})\right.\right\}.
\]
We denote $(f_1^R,f_2^R)$ the solution of \eqref{eq:ODE of norm R} with $R=1/\varepsilon$ and $(f_1,f_2)$ the solution of \eqref{eq:ODE of norm infty}. Then we have
\begin{equation}\label{eq:Ie=ER}
	I(\varepsilon,1)=\pi E_R(f_1^R,f_2^R).
\end{equation}
Lemma \ref{lem:f=s+Cr2+Cr4} and Lemma \ref{lem:f'=Cr3} give that there exist $R_0>0,C>0$ such that for $r>R_0,R>R_0$
\begin{equation}
	\left|f_1^R(r)-\frac{1}{\sqrt{1+g}}\right|\le \frac{C}{r^2},\left|f_2^R(r)-\frac{1}{\sqrt{1+g}} \right|\le \frac{C}{r^2},
\end{equation}
\begin{equation}\label{eq:fj=s+cr2}
	\left|f_1(r)-\frac{1}{\sqrt{1+g}}\right|\le \frac{C}{r^2},\left|f_2(r)-\frac{1}{\sqrt{1+g}} \right|\le \frac{C}{r^2},
\end{equation}
\begin{equation}\label{eq:fj'=Cr3}
	|f'_1(r)|\le \frac{C}{r^{3}},|f_2'(r)|\le \frac{C}{r^{3}}.
\end{equation}
Define 
\begin{equation}
\tilde{f}_j(r)=
	\begin{cases}
	f_j(r),&r\in[0,R-1),\\
	(R-r)f_j(R-1)+(r-R+1),&r\in[R-1,R].
	\end{cases}
\end{equation}
Noting $(f_1^R,f_2^R)$ is a minimizer of $E_R$ in $M_R$, we have
\begin{equation}\label{eq:ERfR<ERf+Cr2}
	E_R(f_1^R,f_2^R)\le E_R(\tilde{f}_1,\tilde{f}_2)\le \int_0^Re_R(f_1,f_2)\ud r+CR^{-2}
\end{equation}
Then define
\begin{equation}
	\tilde{f}_j^R(r)=
	\begin{cases}
	f_j^R(r),&r\in[0,R-1),\\
	f_j^R(r)+(r-R+1)(f_j(R)-f_j^R(R)),&r\in[R-1,R].
	\end{cases}
\end{equation}
Noting that $(f_1,f_2)$ is a  minimizer of $E_R$ in 
\[
	\bar{M}_R:=\{(h_1,h_2):(0,R)\to \R^2|\sqrt{r}h_1',h_1/\sqrt{r},h_2'\in L^2(0,R),h_1(R)=f_1(R),h_2(R)=f_2(R)\},
\]
and $(\tilde{f}_1^R,\tilde{f}_2^R)\in \bar{M}_R$, we have 
\begin{equation}\label{eq:ERf<ERfR+r2}
	E_R(f_1,f_2)\le E_R(\tilde{f}_1^R,\tilde{f}_2^R)\le E_R(f_1^R,f_2^R)+CR^{-2}.
\end{equation}
We have
\begin{equation}
	E_R(f_1,f_2)-\frac{1}{1+g}\log{R}{\sqrt{1+g}}=\int_0^{\sqrt{1+g}}e_R(f_1,f_2)\ud r+\int_{\sqrt{1+g}}^R(e_R(f_1,f_2)-\frac{1}{r(1+g)})\ud r.
\end{equation}
Noting \eqref{eq:fj=s+cr2} and \eqref{eq:fj'=Cr3}, we have
\begin{equation}
	e_R(f_1,f_2)-\frac{1}{r(1+g)}=O(r^{-3}),
\end{equation}
and hence
\begin{equation}\label{eq:ERf-logR=OR2}
	E_R(f_1,f_2)-\frac{1}{1+g}\log{R}{\sqrt{1+g}}=\int_0^{\sqrt{1+g}}e_R(f_1,f_2)\ud r+\int_{\sqrt{1+g}}^{+\infty}(e_R(f_1,f_2)-\frac{1}{r(1+g)})\ud r+O(R^{-2}).
\end{equation}
Combining \eqref{eq:Ie=ER}, \eqref{eq:ERfR<ERf+Cr2}, \eqref{eq:ERf<ERfR+r2} and \eqref{eq:ERf-logR=OR2}, we can prove \eqref{eq:gammag=Ie-loge+e2} with
\begin{equation}
	\gamma_g=\pi\int_0^{\sqrt{1+g}}e_R(f_1,f_2)\ud r+\pi\int_{\sqrt{1+g}}^{+\infty}(e_R(f_1,f_2)-\frac{1}{r(1+g)})\ud r.
\end{equation}
\end{proof}

\begin{Rmk}
\begin{equation}
	I(\varepsilon,R)=\inf_{(u,v)\in A_g(B_R(\bvec{0}))}\int_{B_R(\bvec{0})}e^\varepsilon_0(u,v)\ud\vx
\end{equation}
with 
\[
	A_g(B_R(\bvec{0}))=\left\{(u,v)\in H^1(B_R(\bvec{0}))\times H^1(B_R(\bvec{0}))\left|(u(\vx),v(\vx))=\frac{1}{\sqrt{1+g}}\left(\frac{z}{|z|},1\right)\  \text{for}\ \vx\in \p B_R(\bvec{0})\right.\right\}.
\]
Then wia rescalling, we have $I(\varepsilon,R)=I(\varepsilon/R,1)=I(1,R/\varepsilon)$ and hence
\begin{equation}
	\gamma_g=I(\varepsilon,R)-\frac{\pi}{1+g}\log\frac{R}{\varepsilon\sqrt{1+g}}+O(\varepsilon^2/R^2).
\end{equation}
\end{Rmk}

Here, we give an rough estimate of $\gamma_{1/2}$ via numerical simulation: we use bvp4c to solve the ODE \eqref{eq:ODE of norm e}. We take $\varepsilon=2^j$ with $j=-12,\cdots,-5$ and $\Delta r=\varepsilon/64$. The initial guess of solution is given as 
\begin{equation}
 	f_1(r)=\frac{r+\varepsilon-|r-\varepsilon|}{2\sqrt{1+g}\varepsilon},f_2(r)=\frac{1}{\sqrt{1+g}}.
 \end{equation} 
 We find that when we take $\varepsilon<1/2^{10}$, the error is controlled very well ($<0.001$) and calculated $\gamma_g$ is stable around $0.5377$.

\end{document}